\documentclass[final]{siamltex}

\usepackage{subfigure}
\usepackage{amsmath, amssymb, amsfonts}
\usepackage{bbm}
\usepackage{epsfig}
\usepackage[noline,boxed,longend,algoruled]{algorithm2e}

\newcommand{\Khat}{\widehat{K}}
\newcommand{\Kp}{K^{+}}
\newcommand{\Km}{K^{-}}

\newcommand{\DeltapS}{\Delta_{S}^{+}}
\newcommand{\DeltamS}{\Delta_{S}^{-}}

\newcommand{\KcoupledS}{K_{S}^{0}}

\newcommand{\muodot}{\mu_t^{\odot,N}}
\newcommand{\muoplus}{\mu_t^{\oplus,N}}
\newcommand{\muominus}{\mu_t^{\ominus,N}}

\newcommand{\majT}{\widehat{T}}

\newcommand{\ie}{i.\,e.\ }

\newcommand{\eq}{eq.\ }

\newcommand{\Xminus}{X_t^-}
\newcommand{\Xplus}{X_t^+}

\newcommand{\innerprod}[2]{\left\langle #1 \, , \, #2 \right\rangle}

\newcommand{\wprob}{w.\,p.}

\newcommand{\Var}{\textsc{V}\textrm{ar}}

\newcommand{\varep}{\epsilon}

\renewcommand{\leq}{\leqslant}
\renewcommand{\geq}{\geqslant}

\newcommand{\la}{\lambda}
\newcommand{\al}{\alpha}

\numberwithin{equation}{section} 

\title{Coupling algorithms for calculating sensitivities of Smoluchowski's Coagulation equation}

\author{{Peter~L.~W.~Man\footnotemark[2]\;\,\footnotemark[3]
\and James~R.~Norris\footnotemark[4]
\and Isma\"el~F.~Bailleul\footnotemark[4]
\and Markus~Kraft\footnotemark[5]}}

\footnotetext[2]{
Department of Chemical Engineering and Biotechnology,
University of Cambridge,
New Museums Site,
Pembroke Street,
Cambridge, CB2 3RA,
UK.}

\footnotetext[3]{
A preliminary form of this paper appeared in conference proceedings as the following article: P.~L.~W.~Man, M.~Kraft and J.~R.~Norris, Coupling Algorithms for Calculating Sensitivities of Population Balances, Numerical Analysis and Applied Mathematics, AIP Conference Proceedings, 1048 (2008), pp.~927--930.}

\footnotetext[4]{
Department of Pure Mathematics and Mathematical Statistics,
Centre for Mathematical Sciences,
University of Cambridge,
Wilberforce Road,
Cambridge, CB3 0WB,
UK.}

\footnotetext[5]{
Corresponding author. Department of Chemical Engineering and Biotechnology, University of Cambridge, New Museums Site, Pembroke Street, Cambridge, CB2 3RA, UK (mk306@cam.ac.uk). This research was supported by EPSRC grant EP/E01772X/1.}

\begin{document}

\maketitle

\begin{abstract}
In this paper, two new stochastic algorithms for calculating parametric derivatives of the solution to the Smoluchowski coagulation equation are presented. It is assumed that the coagulation kernel is dependent on these parameters. The new algorithms (called `Single' and `Double') work by coupling two Marcus-Lushnikov processes in such a way as to reduce the difference between their trajectories, thereby significantly reducing the variance of central difference estimators of the parametric derivatives. In the numerical results, the algorithms are shown have have a O(1/N) order of convergence as expected, where N is the initial number of particles. It was also found that the Single and Double algorithms provide much smaller variances. Furthermore, a method for establishing `efficiency' is considered, which takes into account the variances as well as CPU run times, and the `Double' is significantly more `efficient' compared to the `Independent' algorithm in most cases.

\end{abstract}



\pagestyle{myheadings}
\thispagestyle{plain}
\markboth{P.~L.~W.~MAN, J.~R.~NORRIS, I.~F.~BAILLEUL AND M.~KRAFT}{
COUPLING ALGORITHMS FOR CALCULATING SENSITIVITIES}


\bigskip
\smallskip

\section{Introduction}

The simplest of pure coagulation processes puts into play chemical species characterised by a single scalar quantity, say their mass, with values in a discrete set, say the positive integers. The evolution of the process is modelled by a differential equation which gives the time evolution of the concentration $\mu_t^\la(x)$ of particles of mass $x \in \mathbbm{N}$. Given a real-valued function $f$ we shall write $(f,\mu_t^\la)$ for $\sum_{x \in \mathbbm{N}} f(x)\mu_t^\la(x)$. Quantities measured by the experimenter (such as moments) are of this form. Smoluchowski's description of the evolution of $\mu_t^\la$ is
\begin{equation}
\label{SmoEq}
\frac{\textrm{d}}{\textrm{d}t} (f,\mu^\la_t) = \frac{1}{2}\sum_{x,y\geq 1}\bigl\{f(x+y)-f(x)-f(y)\bigr\}\,K_\la(x,y)\,\mu_t^\la(x)\,\mu_t^\la(y).
\end{equation}
The kernel $K_\la(\cdot,\cdot)$ is a symmetric non-negative function which represents the rate at which a pair of particles of masses $x$ and $y$ coagulate to create a particle of mass $x+y$. The term $\{f(x+y)-f(x)-f(y)\bigr\}$ is the change which has occurred in the quantity $(f,\mu_t^\la)$ as a result of this coagulation. The letter $\la$ in the kernel stands for a $d$-dimensional parameter. Our aim in this article is to devise a new numerical scheme for investigating how the solution $\mu_t^\la$ to Smoluchowski equation depends on $\la$. We shall concentrate on the case of a one dimensional parameter since the same analysis applies to the partial derivatives for a multidimensional parameter.

\newpage

There is a large amount of literature concerning the solving of the continuous particle sized version of \eq (\ref{SmoEq}) and its many variations such as particle inception, surface growth, sintering, and fragmentation \cite{Wagner}-\cite{ZhaoZheng09}. However, little is devoted to a systematic method of sensitivity analysis other than merely simulating the physical system in question for various parameter values to measure the change in some quantity as a result of the parameter change. In this paper, we wish to conduct sensitivity analysis by explicitly calculating the parametric derivative of equation~(\ref{SmoEq}).

There are only a few sources which report on this approach \cite{Vik04}-\cite{AsmussenGlynn}. One method involves a weighted particle method, which assigns to each particle a weight with the interpretation of the number of physical particles it represents. This particular method \cite{Vik04,Vik06b} considers a finite difference approach where both particle systems (with different parameters) are simulated together, the only difference being the particles' weights. A potentially very powerful method based on the Lagrangian formalism is considered in \cite{Vik06}. The idea is to consider an adjoint equation which solves for the parametric derivative directly rather than \eq (\ref{SmoEq}). This allows the solving of the derivative for all values of the parameter simultaneously.

\noindent

The aim of this paper is to present two new stochastic algorithms for the calculation of parametric derivatives of \eq (\ref{SmoEq}) with emphasis on variance reduction. These algorithms are based on the simple Marcus--Lushnikov process, but we consider how two such processes with different parameters can be solved simultaneously in order to reduce the estimator variance. These algorithms are presented in section \ref{norris:algorithm}. Their mathematical formalisation is detailed in section \ref{SectionMesaureTheory}, in which we describe how this formalism can be used to justify that these algorithms do indeed provide approximations of the sensitivity. The quality of these approximations is investigated in section \ref{SectionNumerics} where numerical results are analysed.

\section{Central difference estimation of parametric derivatives}
\label{norris:algorithm}

As Marcus-Lushnikov's process is our main ingredient, let us recall first what it is. Dropping the index $\la$, equation \eqref{SmoEq} makes it clear that $\mu_t$ should be seen as a non-negative discrete measure on $\mathbbm{N}$ and $(f,\mu_t) = \sum_{x \in \mathbbm{N}} f(x)\mu_t(x)$ as the integral of $f$ against $\mu_t$. In Marcus-Lushnikov's approach, $\mu_t$ is approximated by a random finite measure of the form\footnote{$\delta_x$ is a Dirac mass at $x \in \mathbbm{N}$.}
\begin{equation*}
\mu_t^N = \frac{1}{N}\sum_{i=1}^n\delta_{x_i(t)}
\end{equation*}
whose dynamics are that of a Markov chain with state space \begin{equation}
\label{eq:simple_state_space}
Q_N:=\left\{ \mu \in \mathcal{M}(\mathbbm{N}) \, \bigg| \, \mu = \frac{1}{N}\sum_{i=1}^n \delta_{x_i} \,,\, x_i \in \mathbbm{N} \,\,\, \forall \,i \, , \, n=1,2,\ldots \, \right\}
\end{equation}
where $\mathcal{M}(\mathbbm{N})$ is the space of all measures on $\mathbbm{N}$. We shall talk of each $x_i(t)$ of $\delta_{x_i(t)}$ as a particle of the system at time $t$. Start from $\mu_0^N = \frac{1}{N}\sum\delta_{x_i}$; associate to each pair $(x_i,x_j)$ of distinct particles an exponential random time $T_{ij}$ with parameter $\frac{K(x_i,x_j)}{N}$, independent of the other exponential times, and set
$$
T:= \min\bigl\{T_{ij}\,;\,i<j\bigr\}.
$$
The process $\mu_t^N$ remains constant on the time interval $[0,T)$ and has a jump at time $T$. If $T=T_{pq}$, set
$$
\mu_T^N = \mu_0^N + \frac{1}{N}\bigl(\delta_{x_p+x_q}-\delta_{x_p}-\delta_{x_q}\bigr);
$$
this operation amounts to removing the particles $x_p$ and $x_q$ from the system and adding the particle $x_p+x_q$. The dynamics then starts afresh.

\smallskip

We shall write $\mu_t^{\la\,;\,N}$ for the Marcus-Lushnikov process corresponding to the kernel $K_\la$. An obvious way of estimating the sensitivity is to approximate it by the (random) ratio $(\mu_t^{\la+\frac{1}{2}\varep\,;\,N}-\mu_t^{\la-\frac{1}{2}\varep\,;\,N})/\varep$. No \textit{a priori} independence or dependence between $\mu_t^{\la+\frac{1}{2}\varep\,;\,N}$ and $\mu_t^{\la-\frac{1}{2}\varep\,;\,N}$ is imposed. We are mainly interested in this article in producing a stochastic approximation of the sensitivity with a low variance. We shall thus try to minimise the variances
$$
\Var\Biggl(\frac{\mu_t^{\la+\frac{1}{2}\varep\,;\,N}(x)-\mu_t^{\la-\frac{1}{2}\varep\,;\,N}(x)}{\varep}\Biggr)
$$
as much as we can for all values of $x$. For that purpose we shall couple the evolution of the two Marcus-Lushnikov processes so as to keep them as close as possible, noting that increasing the covariance between $\mu_t^{\la+\frac{1}{2}\varep\,;\,N}(x)$ and $\mu_t^{\la-\frac{1}{2}\varep\,;\,N}(x)$ decreases the variance. For notational simplicity, we rename $\mu_t^{\la \pm \frac{1}{2}\varep\,;\,N}(x)$ as $\mu_t^{\pm, N}(x)$. It is perfectly possible to describe both trajectories by the $\mathbbm{R}^2$-valued discrete measure $(\mu_t^{+,N},\mu_t^{-,N}) \in Q_N^2$, but to encapsulate the following coupling, we consider the following approach.

\smallskip

\subsection{Coupling}
\label{SectionCoupling}

Here is an example of coupling with framework the unit square of the plane. Denote by $f(x,y)$ any probability density on the square, and consider the problem of minimizing $I:=\int_0^1\int_0^1 |x-y|\, f(x,y)\,\textrm{d}x \, \textrm{d}y$, subject to the condition that the two \emph{marginals} of the probability $f(x,y)\,\textrm{d}x \, \textrm{d}y$ on both the $x$ and $y$ axes are uniform\footnote{That is the marginal $\int_0^1 f(x,y)\textrm{d}y=1$ for each $x \in [0,1]$, and the marginal $\int_0^1 f(x',y)\,\textrm{d}x'=1$ for each $y \in [0,1]$.}. Any measure on the square satisfying this condition is said to realise a \textit{coupling} between the uniform probability on the $x$-segment $[0,1]$ and the uniform probability on the $y$-segment $[0,1]$ (regardless of the above optimisation problem). The probability $\textrm{d}x \, \textrm{d}y$ is such a coupling, but it does not minimise $I$. This minimum is attained for the singular probability on the square with support on the diagonal, and uniform on it, signifying maximum correlation between the $x$ and $y$ axes.

\medskip

Our framework is more complicated than above as the role of $[0,1]$ is now played by the set of particle approximations, i.e. trajectories $\bigl(\{\mu_t^{\pm,N}\}_{t\in [0,t_{\text{end}}]} \bigr)$ (for some final time $t_{\text{end}}$) with values in the set $Q_N$ of finite measures of the form $\frac{1}{N}\sum_i \delta_{x_i}$; but the basic idea is the same. The minimisation of $I$ is replaced by the minimisation of difference between trajectories as seen in the following paragraphs.

\smallskip

Denote by $X^-_t$ and $X^+_t$ the set of particles\footnote{For a measure $\mu = \frac{1}{N} \sum_{i=1}^n \delta_{x_i}$, the set of particles is $\frac{1}{N}(\delta_{x_1},\ldots,\delta_{x_n})$ for $x_i \in \mathbbm{N}$ for all $i$, thus allowing multiple particles with the same size to exist in the set.}
from $\mu_t^{+,N}$ and $\mu_t^{-,N}$, respectively. (The more particles $X^-_t$ and $X^+_t$ have in common the closer the particle systems are, and thus increasing the correlation; the set
\begin{equation}
\label{eq:Xdot_def}
X^{\odot}_t := X^-_t \cap X^+_t
\end{equation}
is made up of those particles in common\footnote{The intersection here is used in the multiset sense - if there are $s_-$ particles of size $x \in \mathbbm{N}$ in $X^-_t$ and $s_+$ particles also of size $x$ in $X^+_t$, then the intersection contains $\min \{s_-,s_+\}$ particles of size $x$. Also, $X_t^{\oplus}$ in eq.~\eqref{eq:Xplusminus_def} then has $\max \{s_+ - s_-,0\}$ particles of size $x$.}, with corresponding measure $\mu_t^{\odot,N} \in Q_N$ being the sum of those particles in $X_t^{\odot,N}$. Similarly, we consider
\begin{equation}
\label{eq:Xplusminus_def}
X^{\oplus}_t := X_t^{+} \backslash X_t^{\odot} \quad, \quad
X^{\ominus}_t := X_t^{-} \backslash X_t^{\odot}
\end{equation}
which are the set of those particles present in $X_t^{+}$ but \textit{not} in $X_t^{-}$ (and the other way round respectively)---we wish to minimise the numbers of particles in these sets. We denote $\mu_t^{\oplus,N}, \mu_t^{\ominus,N} \in Q_N$ as the measures corresponding to $X_t^{\oplus},X_t^{\ominus}$. From this definition, we can recover $\mu_t^{+,N}$ from $\mu_t^{+,N} = \mu_t^{\odot,N} + \mu_t^{\oplus,N}$ and similarly for $\mu_t^{-,N}$. See that the information held in $(X_t^{+},X_t^{-})$ is the same as that held in $(X_t^{\oplus},X_t^{\odot},X_t^{\ominus})$---thus we seek to describe the full stochastic process by the $\mathbbm{R}^3$-valued measure $(\mu_t^{\oplus,N},\mu_t^{\odot,N},\mu_t^{\ominus,N}) \in Q_N^3$ rather than $(\mu_t^{+,N},\mu_t^{-,N}) \in Q_N^2$. Note that the coagulations for $X_t^{\pm}$ particles are governed by the rates determined by $\lambda \pm \frac{1}{2} \epsilon$, and that certain coagulations, such as between a particle in $X_t^{\oplus}$ and a particle in $X_t^{\odot}$, cannot occur in the $X_t^{-}$ set because the particle in $X_t^{\oplus}$ does not exist in $X_t^{-}$. Furthermore, coagulations between a particle in $X_t^{\oplus}$ and a particle in $X_t^{\ominus}$ cannot occur at all.

The first version of our algorithm, called \textbf{Single Coupling Algorithm}, tries to keep the number of particles from $X^{\odot}_t$ as large as possible, imposing that (as much as possible) when two particles, \textit{both} of which are present in $X^-_t$ and $X^+_t$, are chosen to coagulate in one of these systems, they also coagulate in the other. The resulting particle must also be present in both $X^-_t$ and $X^+_t$, thus helping to keep $X^{\odot}_t$ large. Of course, as the coagulation rates in $X^-_t$ and $X^+_t$ differ, we cannot prevent a coagulation event of the above kind from happening in only one of the systems; we can however minimise the rate at which it happens.

\smallskip

The \textbf{Double Coupling Algorithm} is a refinement of the previous one in which we try to make the creation of particles of $X_t^{\oplus}$ and $X_t^{\ominus}$ as rare as possible. In addition to the above coupling, it considers what happens when a particle from $X_t^{\odot}$ coagulates with a particle from $X_t^{\oplus / \ominus}$ (which can only occur in $X_t^{\pm}$ set as mentioned earlier). In such an event, the same particle from $X_t^{\odot}$ can be used in a coagulation event with a particle from $X_t^{\ominus / \oplus}$. Out of the three particles from $X_t^{\oplus}$, $X_t^{\odot}$ and $X_t^{\ominus}$, the $X_t^{\odot}$ particle contributes size to the other two particles, and is itself removed. More details about both couplings are described later.

\smallskip

One ultimately expects that given enough time, the two systems will behave almost independently (\ie{} there will be few particles in $X_t^{\odot}$), but the hope is that the divergence in their trajectories is slow enough over the time span of interest. The simulation of the sensitivity using two independent Marcus-Lushnikov processes will be referred to as the \textbf{Independent Algorithm}; it will be used for comparison with the other algorithms.

\smallskip

\noindent \textbf{Labelling.} The usage of the triple measure $(\mu_t^{\oplus,N},\mu_t^{\odot,N},\mu_t^{\ominus,N}) \in Q_N^3$ captures the similarities and differences between the $\mu_t^{+,N}$ and $\mu_t^{-,N}$ trajectories. Furthermore, all particles are stored one single array, and membership of each particle in one of the particle sets $X_t^{\oplus},X_t^{\odot},X_t^{\ominus}$ is implemented by attaching the particle with a \textbf{label}---these being $\oplus,\odot,\ominus$ respectively. The resulting possible Markov steps are given in Table \ref{List_of_reactions}.

\begin{table}[!h]
\centering
\caption{Particle labels and their meaning}
\begin{tabular}{|l | l|}
\hline
Label & Meaning\\
\hline
$\oplus_x$ & a real particle (of size $x$) present \emph{only} in $X_t^{+}$ \ie{} present in $X_t^{\oplus}$.\\

$\ominus_x$ & a real particle (of size $x$) present \emph{only} in $X_t^{-}$ \ie{} present in $X_t^{\ominus}$.\\

$\odot_x$ & a computational particle (of size $x$) present in $X_t^{\odot}$ \ie{}\\ & a pair of identical real particles, one in $X_t^{+}$ and one in $X_t^{-}$.\\
\hline
\end{tabular}
\label{ParticleLabels}
\end{table}

\begin{table}[!h]
\caption{Possible events described using the labelling notation. 
}
\label{List_of_reactions}
\centering
\renewcommand{\arraystretch}{1}
\begin{tabular}{|p{1.0cm} |l |p{5.5cm}|}
\hline
Type & Event & Explanation\\
\hline\hline
$1 \:
\begin{array}{l}
(a)\\(b)\\(c)
\end{array} $ &
    $
    \odot_{x} + \odot_{y} \rightarrow \left\{
    \begin{array}{l}
    \odot_{x+y}\\
    \ominus_{x} + \ominus_{y} + \oplus_{x+y} \\
    \oplus_{x} + \oplus_{y} + \ominus_{x+y}\\
    \end{array} \right.
    $
     &
     $
     \begin{array}{l}
     \textrm{if occurs in \textbf{both} $\Xminus$ and $\Xplus$}\\
    \textrm{if occurs \textbf{only in} $\Xplus$}\\
    \textrm{if occurs \textbf{only in} $\Xminus$}\\
    \end{array}
    $
    \\
\hline
2(a) & $\oplus_{x} + \odot_{y} + \ominus_{z} \rightarrow \oplus_{x+y} + \ominus_{y+z}$ & See Double Coupling algorithm explanation. Only occurs in Double Coupling algorithm.\\
\hline
2(b) & $\oplus_{x} + \odot_{y} \rightarrow \oplus_{x+y} + \ominus_{y}$ & This represents a coagulation between a pair of particles from $(X_t^{\oplus},X_t^{\odot})$, so the $\oplus$ particle must increase in size and the $\odot$ particle becomes a $\ominus$ (since this particle is no longer in the $\Xplus$ system).\\
\hline
2(c) & $\ominus_{x} + \odot_{y} \rightarrow \ominus_{x+y} + \oplus_{y}$
& Same logic as 2(b)---this reaction can only happen in the $\Xminus$ system.\\
\hline
3(a) & $\oplus_{x} + \oplus_{y} \rightarrow \oplus_{x+y}$ &
Particles present in $X_t^{\oplus}$ coagulate.\\
\hline
3(b) & $\ominus_{x} + \ominus_{y} \rightarrow \ominus_{x+y}$ & As in event type 3(a) except for $X_t^{\ominus}$. \\
\hline
\multicolumn{3}{c}{}\\
\hline
reject & $\oplus_{x} +\ominus_{y}$ & This coagulation cannot occur since each of the particles cannot `see' the other.\\
\hline
\end{tabular}
\end{table}

Note that there is a certain degeneracy in the state space $Q_N^3$---if there exist particles $\frac{1}{N}\delta_{x}$ in \textbf{both} $\mu_t^{\oplus,N}$ and $\mu_t^{\ominus,N}$ then these two particles can be removed and a single particle $\frac{1}{N}\delta_{x}$ added to $\mu_t^{\odot,N}$. Note that this \textbf{cleanup operation} is not a Markov jump but simply a computational enforcement of the definition of $X_t^{\odot}$, and does not affect $X_t^{\pm}$ at all.


\subsection{Single Coupling system}
\label{section:single:coupling}

To be consistent with the above $\pm$ notations, we shall write $K^-$ for the kernel $K_{\la-\frac{1}{2}\varep}$ and $K^+$ for the kernel $K_{\la+\frac{1}{2}\varep}$. Recall that all particles are stored in a single array---we introduce the sets of indices in this array which correspond to particles in $X_t^{\oplus},X_t^{\odot},X_t^{\ominus}$ to be $I(X_t^{\oplus}),I(X_t^{\odot}),I(X_t^{\ominus}) \subseteq \{1,\ldots,n\}$ respectively, where $n$ is the total number of particles\footnote{Thus, $I(X_t^{\oplus}) \cup I(X_t^{\odot}) \cup I(X_t^{\ominus}) = \{1,\ldots,n\}$, and their intersections are empty.}. We then denote $x_i$ to be the size of particle $i$ for $i \in \{1,\ldots,n\}$.

\subsubsection{The idea}

The coupling procedure is implemented using a \textbf{majorant kernel}. This is a symmetric non-negative function $\widehat K(\cdot,\cdot)$ satisfying $\widehat{K}(\cdot,\cdot) \geq K^{\pm}(\cdot,\cdot)$. Run both systems $X_t^-$ and $X_t^+$ at the \textit{same rate}, given by $\widehat K$; a coagulation happening at that rate is called \textit{potential}. If a potential coagulation between particles of sizes $x_i$ and $x_j$ \emph{only} happens in $X_t^{\pm}$, perform it with the respective probabilities
\begin{equation}
\label{reaction_probabilities}
p_{\oplus} = \frac{\Kp(x_i,x_j)}{\Khat(x_i,x_j)}
\quad \textrm{or } \quad
p_{\ominus} = \frac{\Km(x_i,x_j)}{\Khat(x_i,x_j)} \quad,
\end{equation}
otherwise leave the system as it is. This way each system behaves as a Marcus-Lushnikov process with the correct rate.

The coupling itself takes place when the potential coagulation involves a pair of $\odot$-particles (where the coagulation can potentially occur in both $X_t^+$ and $X_t^-$ systems). In this case, the \textit{same} uniform random variable on $(0,1)$ is used to decide whether or not we perform the coagulation event in each system. In other cases the potential coagulation involves only one system. More explicitly, consider only those pairs $(i,j)$ of $\odot$-particles (possibly) involved in the potential coagulation event\footnote{\ie{,} $i,j \in I(X_t^{\odot})$.}.

\begin{figure}[!h]
  \centering
  \includegraphics[width=0.55\textwidth]{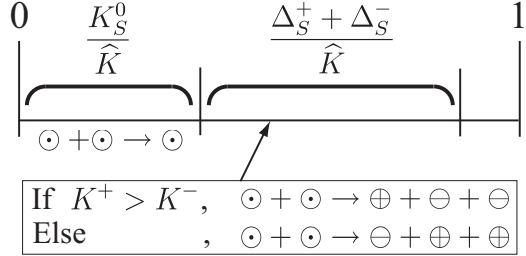}
  \caption{Rate correction for the Single Coupling---generate $\mathbf{U}\sim U(0,1)$ and perform jump event according to the given probabilities.}
  \label{uniform_choosing}
\end{figure}

Set
\begin{romannum}
\item $\KcoupledS(x_i,x_j):=\min\{\Kp(x_i,x_j),\Km(x_i,x_j)\}$ --- rate at which a coagulation of the type $\odot_{x_i} + \odot_{x_j} \rightarrow \odot_{x_i+x_j}$ occurs,
\item $\DeltapS(x_i,x_j):=\max\{\Kp(x_i,x_j)-\Km(x_i,x_j),0\}$ --- rate at which a coagulation of the type $\odot_{x_i} + \odot_{x_j} \rightarrow \ominus_{x_i}+\ominus_{x_j}+\oplus_{x_i+x_j}$ occurs,
\item $\DeltamS(x_i,x_j):=\max\{\Km(x_i,x_j)-\Kp(x_i,x_j),0\}$ --- rate at which a coagulation of the type $\odot_{x_i} + \odot_{x_j} \rightarrow \oplus_{x_i}+\oplus_{x_j}+\ominus_{x_i+x_j}$ occurs.
\end{romannum}
Figure \ref{uniform_choosing} gives a schematic picture of the procedure.

\incmargin{2em}
\Setvlineskip{0.75cm}
\setlength{\algomargin}{1.2em}
\SetKwFor{While}{while}{do}{endwhile}
\SetKwSwitch{Switch}{Case}{Other}{switch}{do}{case}{otherwise}{endswitch}
\begin{algorithm}
\dontprintsemicolon

\textbf{\textsf{For simplicity of exposition, we suppose particles from $X_t^{\oplus},X_t^{\odot},X_t^{\ominus}$ are \emph{all} stored in one single array, which is indexed by $i$ and $j$ below.}}\;

\BlankLine\; \BlankLine\;

\nl Set $t=0$. Set all $N$ initial particles to have $\odot$ labels.\;

\BlankLine\;

\While{$t < t_{\text{end}}$}
{
    \BlankLine\;

    \nl Generate a realisation of the holding time $\Delta t \sim \textrm{Exp}(\widehat{\rho})$ where $\widehat{\rho}$ is specified in \eq \ref{rho_i:def:single} and \eq \ref{rho:def:single}, and set $t \leftarrow t + \Delta t$ .\;

    \BlankLine\; \BlankLine\;

    \textsf{The following step simultaneously chooses the process $k\in\left\{1,2b,2c,3a,3b\right\}$ and the particle pair with indices $(i,j)$ which have the correct labels for the process of type $k$.}\;
    \lnl{algstep:Single:ind_dist} Generate an unordered pair of particles with indices $(i,j)$ for potential coagulation according to the \textbf{index distribution}
    \begin{equation}
    \label{Single_index_dist}
    \frac{\Khat(x_i,x_j)}{2N\widehat{\rho}} \quad,
    \end{equation}
    where particles $i,j$ do not have opposite signs (\ie belong to $I(X_t^{\oplus})$ and $I(X_t^{\ominus})$ respectively, or the other way round). See section \ref{sec:Single_Implementation} for more details.

    \BlankLine\; \BlankLine\;

    \Switch{the value of $k$ chosen}
    {
        \BlankLine\;

        \nl \uCase{$k=1$}
        {
            \textsf{This case represents the Single Coupling part of the algorithm and the following steps (\ref{algstep:Single:Case1_first}-\ref{algstep:Single:Case1_last}) exactly identify with Figure \ref{uniform_choosing}.}\;

            \lnl{algstep:Single:Case1_first} Generate random variable $\textbf{U} \sim U(0,1)$.\;

            \uIf{$0 < \Khat \textbf{U} \leq \KcoupledS$}
            {
                \nl perform event type 1a: $\odot_{x_i} + \odot_{x_j} \rightarrow \odot_{x_i+x_j}$.\;
            }
            \ElseIf{$\KcoupledS < \Khat \textbf{U} \leq \KcoupledS + \DeltapS + \DeltamS$}
            {
                \uIf{$\Kp > \Km$}
                {
                    \nl perform event type 1b: $\odot_{x_i} + \odot_{x_j} \rightarrow \ominus_{x_i} + \ominus_{x_j} + \oplus_{x_i+x_j}$.\;
                }
                \Else
                {
                    \nl perform event type 1c: $\odot_{x_i} + \odot_{x_j} \rightarrow \oplus_{x_i} + \oplus_{x_j} + \ominus_{x_i+x_j}$.\;
                }
            }
            \lnl{algstep:Single:Case1_last} break.\;
        }

        \BlankLine\; \BlankLine\;

        \textsf{For following cases 2b, 2c, 3a and 3b, assume the labels on the particle pair $(i,j)$ match with those as described in steps \ref{algstep:Single:case_2b}-\ref{algstep:Single:case_3b}, otherwise swap the indices $i$ and $j$. Let `\wprob{}' mean `with probability'. Recall also the definitions of $p_{\oplus}$ and $p_{\ominus}$ from eq.~\eqref{reaction_probabilities}.}\;
        \lnl{algstep:Single:case_2b}\lCase{$k=2b$}{ \wprob{~} $p_{\oplus}$, perform $\oplus_{x_i} + \odot_{x_j} \rightarrow \oplus_{x_i+x_j} + \ominus_{x_j}$}. break.\;
        \lnl{algstep:Single:case_2c} \lCase{$k=2c$}{ \wprob{~} $p_{\ominus}$, perform $\ominus_{x_i} + \odot_{x_j} \rightarrow \ominus_{x_i+x_j} + \oplus_{x_j}$}. break.\;
        \lnl{algstep:Single:case_3a}\lCase{$k=3a$}{ \wprob{~} $p_{\ominus}$, perform $\ominus_{x_i} + \ominus_{x_j} \rightarrow \ominus_{x_i+x_j}$}. break.\;
        \lnl{algstep:Single:case_3b}\lCase{$k=3b$}{ \wprob{~} $p_{\oplus}$, perform $\oplus_{x_i} + \oplus_{x_j} \rightarrow \oplus_{x_i+x_j}$}. break.\;
    }

    \BlankLine\; \BlankLine\; \BlankLine\;

    \lnl{algstep:Single:cleanup_step} If a coagulation occurred, for \textbf{each} particle that has just been involved in the coagulation, or newly formed, search for a particle of the same size of the `opposite sign'. If there is such a particle (of size $x$, say), perform a \textit{cleanup operation}: $\ominus_x + \oplus_x \rightarrow \odot_x$.\;

    \BlankLine\;

    \nl \lIf{there is only one particle left in the system}
    {
        STOP.\;
    }

    \BlankLine\;

}

\BlankLine\; \BlankLine\;

\nl STOP.\;

\BlankLine\;

\caption{Single Coupling algorithm}
\label{Single_Coupling_Algorithm_Description}
\end{algorithm}

\subsubsection{The Algorithm}

Recall the different types of coagulation that can happen in the Single Coupling algorithm; they were named $1, 2b, 2c, 3a,$ and $3b$. The total rate of \emph{potential} coagulation is defined as \begin{equation}
\label{rho:def:single} \widehat{\rho}:= \widehat{\rho_1} + \widehat{\rho_{2b}} + \widehat{\rho_{2c}} + \widehat{\rho_{3a}} + \widehat{\rho_{3b}}.
\end{equation}
where the $\widehat\rho_k$ represent the majorant rates at which a \emph{potential} coagulation of type $k \in \{1, 2b, 2c, 3a, 3b\}$ happens:
\begin{subequations}
\label{rho_i:def:single}
\begin{align}
\label{rho_i:def:single_1}
\widehat{\rho_1}&:=\frac{1}{2N}\sum_{\substack{i \neq i'\\i,i' \in I(X_t^{\odot})}}\Khat(x_i,x_{i'}) \quad , &  &\\
\label{rho_i:def:single_2}
\widehat{\rho_{2b}}&:=\frac{1}{N}\sum_{\substack{i,j\\i \in I(X_t^{\odot})\\j \in I(X_t^{\oplus})}}\Khat(x_i,x_j) \quad , &
\widehat{\rho_{2c}}&:=\frac{1}{N}\sum_{\substack{i,k\\i \in I(X_t^{\odot})\\k \in I(X_t^{\ominus})}}\Khat(x_i,x_k)
\quad, \\
\label{rho_i:def:single_3a_3b}
\widehat{\rho_{3a}}&:=\frac{1}{2N}\sum_{\substack{j \neq j'\\j,j' \in I(X_t^{\oplus})}}\Khat(x_j,x_{j'})
\quad , &
\widehat{\rho_{3b}}&:=\frac{1}{2N}\sum_{\substack{k \neq k'\\k,k' \in I(X_t^{\ominus})}}\Khat(x_k,x_{k'})
\quad . &  &
\end{align}
\end{subequations}
Adopting this notation, one can read the details in Algorithm \ref{Single_Coupling_Algorithm_Description} on page \pageref{Single_Coupling_Algorithm_Description}.

\subsubsection{Implementation and Complexity}
\label{sec:Single_Implementation}

The main implementation issue deals with how step \ref{algstep:Single:ind_dist} of Algorithm \ref{Single_Coupling_Algorithm_Description} is performed. We make the assumption that $\widehat{K}$ can be expressed as (for some $A$):
\begin{equation}
\Khat(x_i,x_j) =: \sum_{\al=1}^A f_\al(x_i) g_\al(x_j) \quad;
\end{equation}
such a form of $\widehat K$ enables an easy performance of step \ref{algstep:Single:ind_dist}. As an example of the factorisability condition, the additive kernel $\widehat{K}(x_i,x_j):= \lambda(x_i + x_j)$ can be expressed as $\widehat{K}(x_i,x_j) = \lambda \,.\, x_i + x_j \,.\, \lambda$, implying that $f_1(x)=\lambda$, $g_1(x)=x$, $f_2(x)=x$ and $g_2(x)=\lambda$. The assumption is not so strict---one need only find a \textit{majorant} kernel with this feature. More details on this majorant kernel factorisation can be found in the articles by Eibeck and Wagner \cite{Wagner,Wagner3} and Kraft and coworkers \cite{GoodKraFict,PattersonKraft07,LPDA}.

To see how step \ref{algstep:Single:ind_dist} is performed, first define $C$ as the set of pairs of \emph{distinct} indices (of particles) such that the pair are \emph{not of opposite sign}. Thus eq.~\eqref{Single_index_dist} can be written as:
\begin{subequations}
\begin{align}
\label{PbaSampling}
\frac{\Khat(x_i,x_j)}{2N\widehat{\rho}}
    &= \frac{\sum_{\al} f_\al(x_i) g_\al(x_j)}{\sum_{(p,q) \in C} \sum_{\al'} f_{\al'}(x_{p}) g_{\al'}(x_{q})}\\
    &= \sum_{\al} \left[ \frac{1}{\sum_{(p,q) \in C} \sum_{\al'} f_{\al'}(x_{p}) g_{\al'}(x_{q})} \; \frac{f_\al(x_i) g_\al(x_j)}{1} \right]\\
    \label{PbaSampling_last_equation}
    &= \sum_\al \left[ \frac{\sum_{(p,q) \in C} f_\al(x_p) g_\al(x_q)}{\sum_{(p,q) \in C} \sum_{\al'} f_{\al'}(x_p) g_{\al'}(x_q)}\;  \left( \frac{f_\al(x_i)}{\sum_p f_\al(x_p)}\;\frac{g_{\alpha}(x_j)}{\sum_{q\,;\, (p,q) \in C} g_{\alpha}(x_q)} \right) \right]
\end{align}
\end{subequations}
where $\sum_{(i,j) \in C} f_\al(x_i) g_\al(x_j) = \sum_i f_\al(x_i) \; \sum_{j\,;\, (i,j) \in C} g_{\alpha}(x_j)$ implies the last equality. Thus the user must choose $\alpha$ according to the first fraction of eq.~\eqref{PbaSampling_last_equation} whilst the last two fractions are for the generation of the pair of particles. The advantage of these methods is two-fold: first we can store the values $f_\al(x_i)$ and $g_\al(x_i)$ in `binary' tree structures which also stores their sums (over $i$). This allows efficient generation from the respective distributions
\begin{equation*}
\frac{f_\al(x_i)}{\sum_{p}f_\al(x_p)}
\quad \textrm{and} \quad
\frac{g_\al(x_i)}{\sum_{q}g_\al(x_q)} \quad.
\end{equation*}
Updating the values in this data structure is efficient. If the number of stochastic particles in a binary tree is $n$, the complexity of updating and generating particle operations take $O(\log n)$ steps. Furthermore, the generation of a particle pair is simple---the factorisation allows one to generate each particle in the pair separately meaning that generation of the pair of particles is $O(\log n)$ rather than $O(n^2)$.

In Step \ref{algstep:Single:cleanup_step} of Algorithm \ref{Single_Coupling_Algorithm_Description}, a search of particles for cleanup is required for each iteration. This can be achieved by maintaining linked lists of information about where particles of certain size and label can be found on the particle ensemble list. In short, the Single Coupling algorithm may be faster than the `Independent' algorithm since we need to simulate for one particle ensemble rather than two. On the other hand, the cleanup procedure in the Single Coupling requires extra storage of information, and computational time to update this information.

The Single Coupling algorithm is good for the initial prevention of creation of $\oplus$ and $\ominus$ particles, however, as the $\odot$ particle numbers decrease over time, the Single Coupling should become less effective. This motivates the Double Coupling procedure which is designed to minimise the rate of creation of $\oplus$ and $\ominus$ particles for later times.


\subsection{Double Coupling system}

\begin{figure}[!h]
\centering
\includegraphics[width=0.7\textwidth]{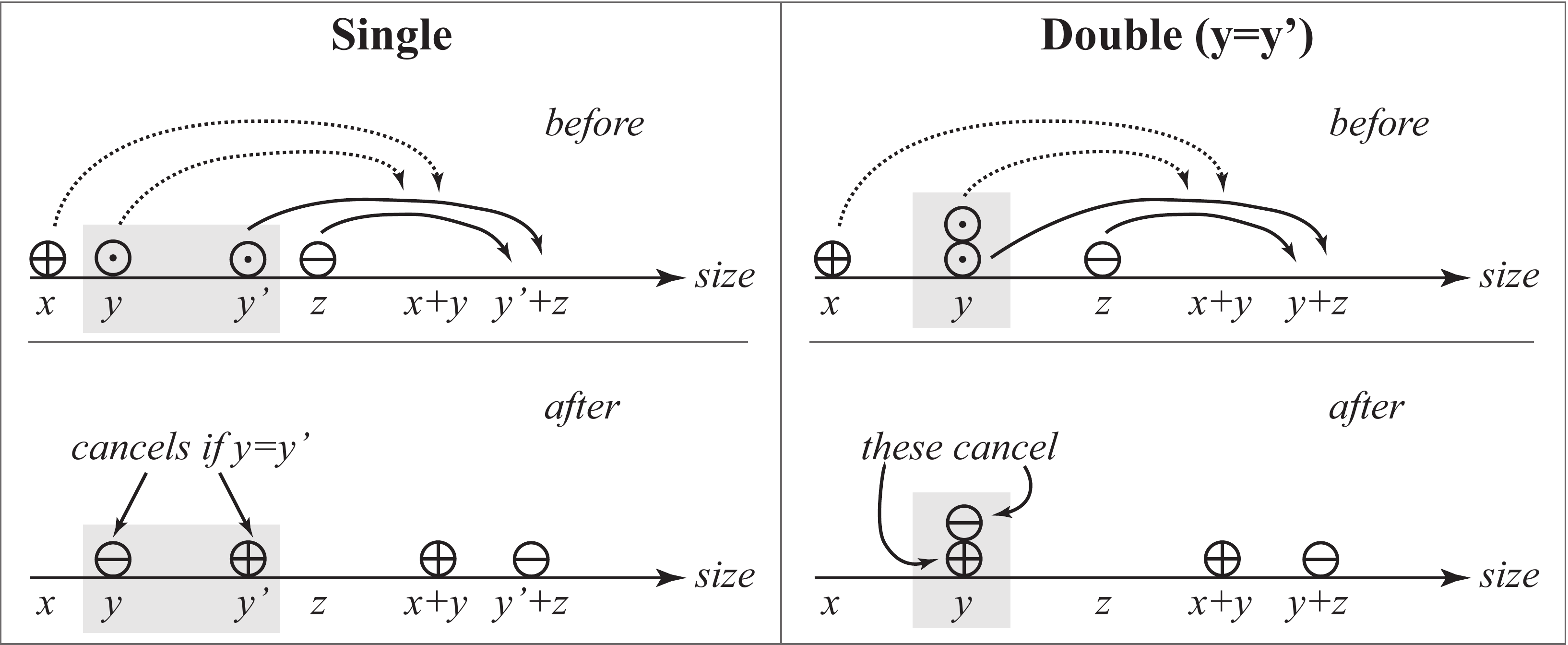}\\
\caption{Pictorial explanation of the Double Coupling algorithm.}
\label{Double_explanation}
\end{figure}

The aim of the Double Coupling algorithm is to try to minimise the rate at which particles of type $\ominus$ or $\oplus$ are created; it was briefly described in section \ref{SectionCoupling}. Figure \ref{Double_explanation} presents a pictorial illustration of this coupling.

\noindent More formally, we
\begin{romannum}
\item choose a $\odot$ particle as the common particle for the $\ominus + \odot$ and $\oplus + \odot$ coagulations. This is done at the maximum potential rate at which the two reactions can occur simultaneously (for a common $\odot$ particle $i \in I(X_t^{\odot})$)
    \begin{displaymath}
    \max\left\{\sum_{j' \in I(X_t^{\oplus})}\Khat(x_{j'},x_i),\sum_{k' \in I(X_t^{\ominus})}\Khat(x_{k'},x_i)\right\},
    \end{displaymath}
\item choose a $\ominus$ particle $k \in I(X_t^{\ominus})$ (for a $\ominus + \odot$ coagulation) \emph{and} a $\oplus$ particle $j \in I(X_t^{\oplus})$(for a potential $\oplus + \odot + \ominus$ coagulation) with respective distributions
    \begin{eqnarray*}
    \frac{\Khat(x_k,x_i)}{\sum_{k' \in I(X_t^{\ominus})}\Khat(x_{k'},x_i)}
    & \quad
    \textrm{and}
    & \quad
    \frac{\Khat(x_j,x_i)}{\sum_{j' \in I(X_t^{\oplus})}\Khat(x_{j'},x_i)}
    \end{eqnarray*}
\item rejection steps are performed to correct the rates according to whether the coagulation happens in each of the $X_t^{+}$ and $X_t^{-}$ systems.
\end{romannum}

\smallskip


\incmargin{2em}
\Setvlineskip{0.75cm}
\setlength{\algomargin}{1.2em}
\SetKwFor{While}{while}{do}{endwhile}
\SetKwSwitch{Switch}{Case}{Other}{switch}{do}{case}{otherwise}{endswitch}
\begin{algorithm}
\dontprintsemicolon

\textbf{\textsf{The Double Coupling algorithm is almost identical to the Single Coupling algorithm, but is modified by \textit{replacing} cases $k=2b$ and $k=2c$ (steps \ref{algstep:Single:case_2b} and \ref{algstep:Single:case_2c}) in Algorithm \ref{Single_Coupling_Algorithm_Description} with a new combined case $k=2$ containing the following steps (ignore any particle pair already chosen)}}.\;

\BlankLine\; \BlankLine\; \BlankLine\;

\nl Choose a $\odot$ particle $i \in I(X_t^{\odot})$ with the distribution
\begin{equation}
\label{distribution_for_neutral_particle}
\frac{\majT_N(+,i)+\majT_N(-,i)}{\sum_{i' \in I(X_t^{\odot})}[\majT_N(+,i')+\majT_N(-,i')]} = \frac{\majT_N(+,i)+\majT_N(-,i)}{N \widehat{\rho_2}} \quad.
\end{equation}\;

\BlankLine\; \BlankLine\; \BlankLine\;

\textsf{\textbf{The first rejection step}---the following step is there purely to transform the total potential rate of process 2 into $\sum_{i' \in I(X_t^{\odot})}\majT_N{(\vee,i')}$ from $\widehat{\rho_2}$.}\;
\nl With probability
\begin{equation*}
\frac{\majT_N{(\vee,i)}}{\majT_N(+,i)+\majT_N(-,i)}
\end{equation*}
we continue, else reject by going to Step \ref{double_reject}.\;

\BlankLine\; \BlankLine\;

\nl Choose a $(\oplus,\ominus)$ particle pair $(j,k)$ with $j \in I(X_t^{\oplus})$ and $k \in I(X_t^{\ominus})$ according to the respective distributions
\begin{equation}
\label{distribution_for_both_particle}
\frac{\Khat(x_j,x_i)}{\majT_N(+,i)}
    \quad \textrm{and} \quad
\frac{\Khat(x_k,x_i)}{\majT_N(-,i)} \quad.
\end{equation}\;

\BlankLine\; \BlankLine\; \BlankLine\;

\textsf{We now have generated a triplet of particles $(i,j,k) \in I(X_t^{\odot}) \times I(X_t^{\oplus}) \times I(X_t^{\ominus})$. Define the probabilities of the $\oplus +\odot$ and $\ominus +\odot$ coagulations occurring in $X_t^{+},X_t^{-}$ respectively as:
\begin{equation}
p_{\oplus+\odot} := \frac{\majT_N(+,i)}{\majT_N(\vee,i)} \frac{\Kp(x_j,x_i)}{\Khat(x_j,x_i)}
\quad \textrm{and} \quad
p_{\ominus+\odot} := \frac{\majT_N(-,i)}{\majT_N(\vee,i)} \frac{\Km(x_k,x_i)}{\Khat(x_k,x_i)} \quad.
\end{equation}
\textbf{The second rejection} (steps \ref{gen_unif}-\ref{dot_minus})---this occurs in an almost identical fashion to Figure \ref{uniform_choosing}, just with different rates.}\;

\lnl{gen_unif} Generate random variable $\textbf{U} \sim U(0,1)$.\;

\uIf{$\textbf{U} < \min \{p_{\oplus + \odot}, p_{\ominus+\odot} \}$}{
    \nl perform event type 2a: $\oplus_x + \odot_y + \ominus_z \rightarrow \oplus_{x+y} + \ominus_{y+z}$.\;
}
\ElseIf{$\min \{p_{\oplus + \odot}, p_{\ominus+\odot} \} \leq \textbf{U} < \max \{p_{\oplus + \odot}, p_{\ominus+\odot} \}$}
{
    \uIf{$p_{\oplus + \odot} > p_{\ominus + \odot}$}
    {
        \nl perform event type 2b: $\oplus_x + \odot_y \rightarrow \oplus_{x+y} + \ominus_{y}$.\;
    }
    \Else
    {
        \lnl{dot_minus} perform event type 2c: $\ominus_{z} + \odot_{y} \rightarrow \ominus_{y+z} + \oplus_{y}$.\;
    }
}

\BlankLine\; \BlankLine\; \BlankLine\;

Go to Step \ref{algstep:Single:cleanup_step} of the Single Coupling algorithm (Algorithm \ref{Single_Coupling_Algorithm_Description}) \lnl{double_reject}.\;

\caption{Double Coupling algorithm - only the part which differs from Algorithm \ref{Single_Coupling_Algorithm_Description}.}

\label{Double_Coupling_Algorithm_Description}
\end{algorithm}

\subsubsection{The algorithm}

The algorithm is the same as the Single Coupling version, except that we merge processes 2b and 2c into a new process 2 whose majorant rate is $\widehat\rho_2:=\widehat\rho_{2b}+\widehat\rho_{2c}$. In describing the algorithm, and given a \emph{particular} $\odot$ particle $i \in I(X_t^{\odot})$, we write $\majT_N(+,i)$ for $\sum_{j' \in I(X_t^{\oplus})}\Khat(x_{j'},x_i)$, and  $\majT_N(-,i)$ for $\sum_{k' \in I(X_t^{\ominus})}\Khat(x_{k'},x_i)$; the maximum of these two quantities is denoted by $\majT_N(\vee,i)$. See Algorithm \ref{Double_Coupling_Algorithm_Description} on page \pageref{Double_Coupling_Algorithm_Description} for the description of the Double Coupling algorithm. In the next paragraph we directly check that the algorithm produces coagulations with the correct rates.

\smallskip

\paragraph{Double Coupling algorithm rates} In the specification of the state space earlier, we recall that $\muoplus,\muodot,\muominus$ are the empirical measures for the $\odot,\oplus,\ominus$ particles at time $t$ respectively. Given a particular $\odot$ particle of mass $y$, the total majorant \emph{potential} rate at which this particle reacts with any $\oplus$ particle is equal to
\begin{equation}
\majT_N(+,y) := \displaystyle \sum_{x \geq 1} \Khat(x,y) \, \muoplus(x) \quad,
\label{majT_measure_def}
\end{equation}
and similarly for $\majT_N(-,y)$ and $\majT_N(\vee,y)$, so that these are analogous quantities to $\majT_N(+,i),\majT_N(-,i),\majT_N(\vee,i)$ used in Algorithm \ref{Double_Coupling_Algorithm_Description}. Note that $\majT_N(\pm,y)$  are functionals depending on $\mu_t^{\oplus / \ominus , N}$ respectively. Also, the majorant rate at which a coagulation event of the form $\oplus_x + \odot_y + \ominus_z \rightarrow \oplus_{x+y} + \ominus_{y+z}$ occurs is\footnote{We use the index $D$ for ``Double''; this distinguishes the quantities to be introduced from the similar ones introduced above for the Single Coupling algorithm.}
\begin{multline*}
\lefteqn{K^{0,N}_D(x,y,z) := }\\
    \underbrace{\widehat{\rho_2}}_{\substack{\textrm{Total rate}\\ \textrm{of process 2}}} \cdot
    \underbrace{\frac{\majT_N(+,y)+\majT_N(-,y)}{\displaystyle \sum_{y \geq 1} \left[ \majT_N(+,y)+\majT_N(-,y) \right] \, \muodot(y)}}_{\textrm{Choose
    $\odot$ particle}} \cdot
    \underbrace{\frac{\majT_N(\vee,y)}{\majT_N(+,y)+\majT_N(-,y)}}_{\textrm{1st Rejection step}} \cdot\\
    \underbrace{\frac{\Khat(x,y)}{\majT_N(+,y)
    }}_{\textrm{Choose a $\oplus_{x}$}} \cdot
    \underbrace{\frac{\Khat(z,y)}{\majT_N(-,y)
    }}_{\textrm{Choose a $\ominus_{z}$}} \cdot
    \underbrace{\min
        \left\{ \frac{\Km(z,y)}{\Khat(z,y)} \, \frac{\majT_N(-,y) }{\majT_N(\vee,y)} ,
        \frac{\Kp(x,y)}{\Khat(x,y)} \, \frac{\majT_N(+,y)}{\majT_N(\vee,y)}
        \right\}
    }_{\textrm{Probability of rejecting neither
    coagulation}}
\end{multline*}
which simplifies to
\begin{equation}
\label{min_expression}
K^{0,N}_D(x,y,z) = \min \left\{
    r^{-}(x,y,z),
    r^{+}(x,y,z)
\right\},
\end{equation}
where
\begin{equation}
r^{-}(x,y,z) := \frac{\Khat(x,y)
    \Km(y,z)}{\majT_N(+,y)}
\quad ,\quad
r^{+}(x,y,z) := \frac{\Khat(y,z)
    \Kp(x,y)}{\majT_N(-,y)} \quad.
\label{joint_double_rates_def}
\end{equation}
Similarly (and dropping the $(x,y,z)$ for convenience), the rate at which \textbf{only} the $\Xminus$ reaction occurs is
\begin{subequations}
\label{Double_non_coupling_rates_def}
\begin{align}
\Delta^{-,N}_D
    &:= \bigl[
            \max \left\{
                 r^{-},r^{+}
            \right\}
            -
            \min \left\{
                r^{-},r^{+}
            \right\}
         \bigr]
        \mathbbm{1}_{r^{+} < r^{-}} \nonumber \\
    &= \max \{ r^{-} - r^{+} , 0 \} \label{Double_non_coupling_rates_def:1}, \\
\intertext{and the rate at which only the $\Xplus$ reaction occurs is}
\Delta^{+,N}_D &:= \max \{ r^{+} - r^{-} , 0 \}
\label{Double_non_coupling_rates_def:2}.
\end{align}
\end{subequations}
As verification for the above rate expressions, we note that the rate at which a pair of particles $(\oplus_x,\odot_y)$ coagulates (in $X_t^+$) is
\begin{align*}
\displaystyle{\sum_{z\geq 1}}\,\bigl(K^{0,N}_D(x,y,z) + \Delta^{+,N}_D(x,y,z)\bigr) \, \muominus(z)
& = \displaystyle{\sum_{z \geq 1}}\;r^{+}(x,y,z) \, \muominus(z) \nonumber \\
& = \Kp(x,y) \, \displaystyle{\sum_{z \geq 1}}\;\frac{\Khat(y,z)}{\majT(-,y)} \, \muominus(z) \nonumber \\
& = \Kp(x,y).
\end{align*}
A similar computation is made to check that $(\odot_y,\ominus_z)$ coagulate in $X_t^-$ at rate $\Km(y,z)$.

\subsubsection{Implementation and Complexity}

Looking at Step \ref{algstep:Single:ind_dist} of the Single Coupling algorithm where the particle pair is chosen, and simultaneously the process $k$, we note that the combined process $\widehat{\rho_2}$ for the Double Coupling is chosen by choosing either a $(\odot,\ominus)$ or a $(\odot,\oplus)$ particle pair. Either way, a $\odot$ particle is automatically chosen with the correct distribution in equation \ref{distribution_for_neutral_particle} and one of $\oplus$ and $\ominus$ is also automatically chosen with the correct distribution specified in equation \ref{distribution_for_both_particle} respectively. This only leaves the remaining particle left to be chosen.

The complexity of this algorithm should be similar to that of the Single Coupling, except that the combined process $\widehat{\rho_2}$ requires slightly more work than in the Single Coupling. However, the Double Coupling hopefully reduces the number of $\oplus$ and $\ominus$ and therefore would reduce the total rate of reactions. Consequently, there might be slightly fewer coagulation events in total.

\subsubsection{Limit coupled processes}
\label{SectionMesaureTheory}
Recall the definitions of the measures $\mu_t^{\oplus,N}$, $\mu_t^{\odot,N}$, $\mu_t^{\ominus,N}$ as given in section \ref{SectionCoupling}. In the same way as one can prove that the Marcus-Lushnikov process converges to the solution of Smoluchowski equation (when it is unique)\footnote{See for instance the article \cite{Nor99} of J. Norris or \cite{Jeon} of I. Jeon.}, it is reasonable to propose a similar result for the triple of \textit{stochastic} processes $\bigl(\mu^{\oplus,N}_\cdot, \mu^{\odot,N}_\cdot, \mu^{\ominus,N}_\cdot\bigr)$. The limiting object $\bigl(\mu^{\oplus}_\cdot, \mu^{\odot}_\cdot, \mu^{\ominus}_\cdot\bigr)$ is a \textit{deterministic} non-negative measure-valued path. Given three bounded functions $f,g,h$, it satisfies the system\footnote{The rates $K_D^0$, $\Delta_D^{\pm}$ and $\widehat{T}(\pm,y)$ are the analogues to $K_D^{0,N}$, $\Delta_D^{\pm,N}$ and $\majT_N(\pm,y)$ but with dependence on the measures $\mu_t^{\oplus}, \mu_t^{\ominus}$ rather than on $\mu_t^{\oplus,N}, \mu_t^{\ominus,N}$. See eqs.~\eqref{min_expression} to \eqref{Double_non_coupling_rates_def} for the expressions for $K_D^{0,N}$, $\Delta_D^{\pm,N}$ and $\majT_N(\pm,y)$.}
\begin{subequations}
\label{Double_weak_eq}
\begin{align}
\label{Double_weak_eq:1}
\frac{\textrm{d}}{\textrm{d}t} \bigl(f,\mu^\odot_t\bigr)
    &= \frac{1}{2}\sum_{x,y\geq 1} \left[ f(x+y)-f(x)-f(y) \right] \, \KcoupledS(x,y) \, \mu^\odot_t(x)\,\mu^\odot_t(y) \nonumber\\
    &- \sum_{x,y\geq 1} f(x)\, \left[\DeltapS(x,y)+\DeltamS(x,y)\right]\,\mu^\odot_t(x)\,\mu^\odot_t(y) \nonumber\\
    &- \sum_{x,y,z\geq 1} f(y)\, K^0_D(x,y,z)\,\mu^\oplus_t(x)\,\mu^\odot_t(y)\,\mu^\ominus_t(z) \nonumber\\
    &- \sum_{x,y,z\geq 1} f(y)\, \left[ \Delta_D^+(x,y,z) + \Delta_D^-(x,y,z) \right] \,\mu^\oplus_t(x)\,\mu^\odot_t(y)\,\mu^\ominus_t(z)\\
\intertext{and}
\label{Double_weak_eq:2}
\frac{\textrm{d}}{\textrm{d}t}\bigl(g,\mu^\oplus_t\bigr)
    &= \frac{1}{2}\sum_{x,y\geq 1} g(x+y)\, \DeltapS(x,y)\,\mu^\odot_t(x)\,\mu^\odot_t(y) + \sum_{x,y\geq 1} g(x)\, \DeltamS(x,y) \mu^\odot_t(x)\mu^\odot_t(y) \nonumber\\
    &+ \sum_{x,y,z\geq 1} \left[ g(x+y)-g(x) \right] \, \bigl(K^0_D+\Delta^+_D\bigr)(x,y,z)\,\mu^\oplus_t(x)\,\mu^\odot_t(y)\,\mu^\ominus_t(z) \nonumber\\
    &+ \sum_{x,y,z\geq 1} g(y) \, \Delta_D^-(x,y,z)\,\mu^\oplus_t(x)\,\mu^\odot_t(y)\,\mu^\ominus_t(z) \nonumber\\
    &+\frac{1}{2}\sum_{x,y\geq 1} \left[ g(x+y)-g(x)-g(y) \right]\, \Kp(x,y)\,\mu^\oplus_t(x)\mu^\oplus_t(y) \quad .
\end{align}
\end{subequations}
A similar equation to eq.~\eqref{Double_weak_eq:2} holds for $\frac{\textrm{d}}{\textrm{d}t}\bigl(h,\mu^\ominus_t\bigr)$. The reader will get a clear insight on the reason why these equations appear by seeing the generator of the discrete measure valued Markov chain $\bigl(\mu_t^{\oplus,N},\mu_t^{\odot,N},\mu_t^{\ominus,N}\bigr)$ (this will be shown in the next subsection \ref{sec:Generator}). Recall that
\begin{equation*}
\mu_t^{+,N} = \mu_t^{\odot,N} + \mu_t^{\oplus,N}, \qquad
\mu_t^{-,N} = \mu_t^{\odot,N} + \mu_t^{\ominus,N},
\end{equation*}
and so it is easily shown that for any bounded functions $f$ and $g$
\begin{equation*}
\frac{\textrm{d}}{\textrm{d}t} (f,\mu_t^+) = \frac{1}{2}  \sum_{x,y\geq 1}\left[ f(x+y)-f(x)-f(y) \right] \, \Kp(x,y) \, \mu_t^+(x) \, \mu_t^+(y)
\end{equation*}
and
\begin{equation*}
\frac{\textrm{d}}{\textrm{d}t} (g,\mu_t^-) = \frac{1}{2} \sum_{x,y\geq 1} \left[ g(x+y)-g(x)-g(y) \right] \, \Km(x,y) \, \mu_t^-(x) \, \mu_t^-(y) \quad.
\end{equation*}
This is in accordance with the fact that $\mu_t^{\pm,N}$ are Marcus-Lushnikov processes with rates $K^{+/-}$, therefore their limits are solutions to Smoluchowski equation with the corresponding rate (under certain conditions). This implies that their difference converges to the difference of the two solutions, this being true independent of the coupling.

\subsubsection{Generator}
\label{sec:Generator}


This section gives a description of the generator of the Markov chain corresponding to the Double Coupling Algorithm. The stochastic jumps are described by the following elementary operations on measures corresponding to the jumps indicated in Figure \ref{List_of_reactions}. We adopt the notations $\mu$ for a generic element of $Q_N^3$ and $x,y,y',z$ for integer masses.
\begin{equation*}
\begin{array}{rcc lllll l}
J_{1a}^N(\mu,y,y') &=& \mu + \frac{1}{N}(& 0 &,& \delta_{y+y'} - \delta_{y} - \delta_{y'} &,& 0 &)\\
J_{1b}^N(\mu,y,y') &=& \mu + \frac{1}{N}(& \delta_{y+y'} &,& -\delta_{y} - \delta_{y'} &,& \delta_{y} + \delta_{y'} &)\\
J_{1c}^N(\mu,y,y') &=& \mu + \frac{1}{N}(& \delta_y + \delta_{y'} &,& -\delta_y - \delta_{y'} &,& \delta_{y+y'} &)\\\\
J_{2a}^N(\mu,x,y,z) &=& \mu + \frac{1}{N}(& -\delta_{x} + \delta_{x+y} &,& -\delta_y &,& -\delta_{z} + \delta_{y+z} &)\\
J_{2b}^N(\mu,x,y,z) &=& \mu + \frac{1}{N}(& -\delta_{x} + \delta_{x+y} &,& -\delta_y &,& \delta_{y} &)\\
J_{2c}^N(\mu,x,y,z) &=& \mu + \frac{1}{N}(& \delta_{y} &,& -\delta_{y} &,& -\delta_{z} + \delta_{y+z} &)\\\\
J_{3a}^N(\mu,x,y) &=& \mu + \frac{1}{N}(& 0 &,& 0 &,& \delta_{x+y} - \delta_{x} - \delta_y &)\\
J_{3b}^N(\mu,x,y) &=& \mu + \frac{1}{N}(& \delta_{x+y} - \delta_{x} - \delta_y &,& 0 &,& 0 &)
\end{array}
\end{equation*}
The introduction of the following notation clarifies the description of the generator of the Markov chain corresponding to the double coupling algorithm. For any $\gamma \in Q_N$, define the rescaled counting measure $\widetilde{\gamma} \in Q_N$ on ordered pairs of masses of distinct particles as
\begin{equation*}
\widetilde{\gamma}(A \times A') := \gamma(A)\gamma(A') - \frac{1}{N}\gamma(A \cap A') \quad, \quad A,A' \subset \mathbbm{N}.
\end{equation*}
Set also for any $\phi = (\phi_1,\phi_2,\phi_3)$ and $\mu=(\mu_1,\mu_2,\mu_3) \in Q_N^3$
\begin{equation*}
\left\langle \phi,\mu \right\rangle := \left( \langle \phi_1,\mu_1\rangle , \langle \phi_2,\mu_2\rangle,  \langle \phi_3,\mu_3\rangle \right).
\end{equation*}
For any measure $\mu := (\mu^{\oplus},\mu^{\odot},\mu^{\oplus}) \in Q_N^3$, with the corresponding \emph{sets} of particles being $(X^{\oplus},X^{\odot},X^{\ominus})$, we have

\begin{align*}
\innerprod{(\phi_1,\phi_2,\phi_3)}{\mathcal{G}_t^N(\mu)} :=
\left( \innerprod{\phi_1}{\mathcal{G}_t^{N,\oplus}(\mu)} \, , \, \innerprod{\phi_2}{\mathcal{G}_t^{N,\odot}(\mu)} \, , \,
\innerprod{\phi_3}{\mathcal{G}_t^{N,\ominus}(\mu)} \right)
\end{align*}
where
\begin{subequations}
\label{eq:Generator1_explicit}
\begin{align}
\label{eq:Generator1_explicit_plus}
\lefteqn{\innerprod{\phi_1}{\mathcal{G}_t^{N,\oplus}(\mu)}} & \nonumber\\
    &= \frac{1}{2}\sum_{y,y'\geq 1} \phi_1(y+y')\, \DeltapS(y,y)\,\widetilde{\mu}^\odot(y,y') + \sum_{y,y'\geq 1} \phi_1(y)\, \DeltamS(y,y') \widetilde{\mu}^\odot(y,y') \nonumber\\
    &+ \sum_{x,y,z\geq 1} \left[ \phi_1(x+y)-\phi_1(x) \right] \, \bigl(K^0_D+\Delta^+_D\bigr)(x,y,z)\,\mu^\oplus(x)\,\mu^\odot(y)\,\mu^\ominus(z) \nonumber\\
    &+ \sum_{x,y,z\geq 1} \phi_1(y) \, \Delta_D^-(x,y,z)\,\mu^\oplus(x)\,\mu^\odot(y)\,\mu^\ominus(z) \nonumber\\
    &+\frac{1}{2}\sum_{x,x'\geq 1} \left[ \phi_1(x+x')-\phi_1(x)-\phi_1(x') \right]\, \Kp(x,x')\,\widetilde{\mu}^\oplus(x,x') \\
\intertext{and}
\label{eq:Generator1_explicit_dot}
\lefteqn{\innerprod{\phi_2}{\mathcal{G}_t^{N,\odot}(\mu)}} & \nonumber\\
    &= \frac{1}{2}\sum_{y,y'\geq 1} \left[ \phi_2(y+y')-\phi_2(y)-\phi_2(y') \right] \, \KcoupledS(y,y') \, \widetilde{\mu}^{\odot}(y,y') \nonumber\\
    &- \sum_{y,y'\geq 1} \phi_2(y)\, \left[\DeltapS(y,y')+\DeltamS(y,y')\right]\,\widetilde{\mu}^{\odot}(y,y') \nonumber\\
    &- \sum_{x,y,z\geq 1} \phi_2(y)\, K^0_D(x,y,z)\,\mu^\oplus(x)\,\mu^\odot(y)\,\mu^\ominus(z) \nonumber\\
    &- \sum_{x,y,z\geq 1} \phi_2(y)\, \left[ \Delta_D^+(x,y,z) + \Delta_D^-(x,y,z) \right] \,\mu^\oplus(x)\,\mu^\odot(y)\,\mu^\ominus(z) \quad,
\end{align}
\end{subequations}
with $\innerprod{\phi_3}{\mathcal{G}_t^{N,\ominus}(\mu)}$ being defined analogously to eq.~\eqref{eq:Generator1_explicit_plus}.

\bigskip

\section{Numerical Results}
\label{SectionNumerics}

The results presented consider two kernels: the additive kernel $K(x,y)=\la(x+y)$ and a kernel that is used in modelling soot formation in a free molecular regime (thus we shall call it the `Soot Kernel')\footnote{This kernel is studied extensively in \cite{GoodKraFict} and used in \cite{PattersonSinghBalthasar,NealMarkusBalthasarWongFrenklachMitchell,CelnikPattersonKraftWagner}.}
\begin{equation*}
K(x,y) = \left(\frac{1}{x}+\frac{1}{y}\right)^{\frac{1}{2}} \left(x^{\frac{1}{\lambda}} + y^{\frac{1}{\lambda}}\right)^2.
\end{equation*}
The reference value of $\la$ for the additive kernel will be $1$ and for the soot kernel $2.1$. We shall always take as initial condition for the Marcus-Lushnikov process $N$ particles with mass equal to $1$. Throughout this section we shall denote by $N$ the initial number of particles in each system (which is the same), by $\la$ the above reference value of the parameter, whose perturbation will be denoted by $\epsilon$ (\ie the $X^\pm$ systems are governed by the parameter values $\la \pm \frac{1}{2}\epsilon$), and by $L$ the number of simulations with the same initial conditions. The remaining notation is given below.
\begin{romannum}
\item $t$ = time of evolution of the particle system
\item $t_{\textrm{run}}$ = time taken to run the algorithms (CPU time).
\item The estimate of $\frac{\partial}{\partial \lambda} (f,\mu_t^{\lambda})$ given by the $l^{\textrm{th}}$ simulation is denoted by $F_l^{(\lambda)}$(\footnote{Note that taking $f(y)={\bf 1}_{y=x}$ gives $ \frac{\partial}{\partial \lambda}\mu_t^{\lambda}(x)$.}), where $f$ is a suitable test function.
\item The estimate of $\frac{\partial}{\partial \lambda} (f,\mu_t^{\lambda})$ given by $L$ simulations is denoted by $\overline{F}^{(\lambda)}$; it is equal to $\frac{1}{L} \sum_{l=1}^L F_l^{(\lambda)}$.
\end{romannum}


\subsection{Some initial plots}

\begin{figure}[!h]
\begin{center}
    \subfigure[][$t=0.5$]
    {
        \resizebox{!}{0.35\linewidth}{\includegraphics{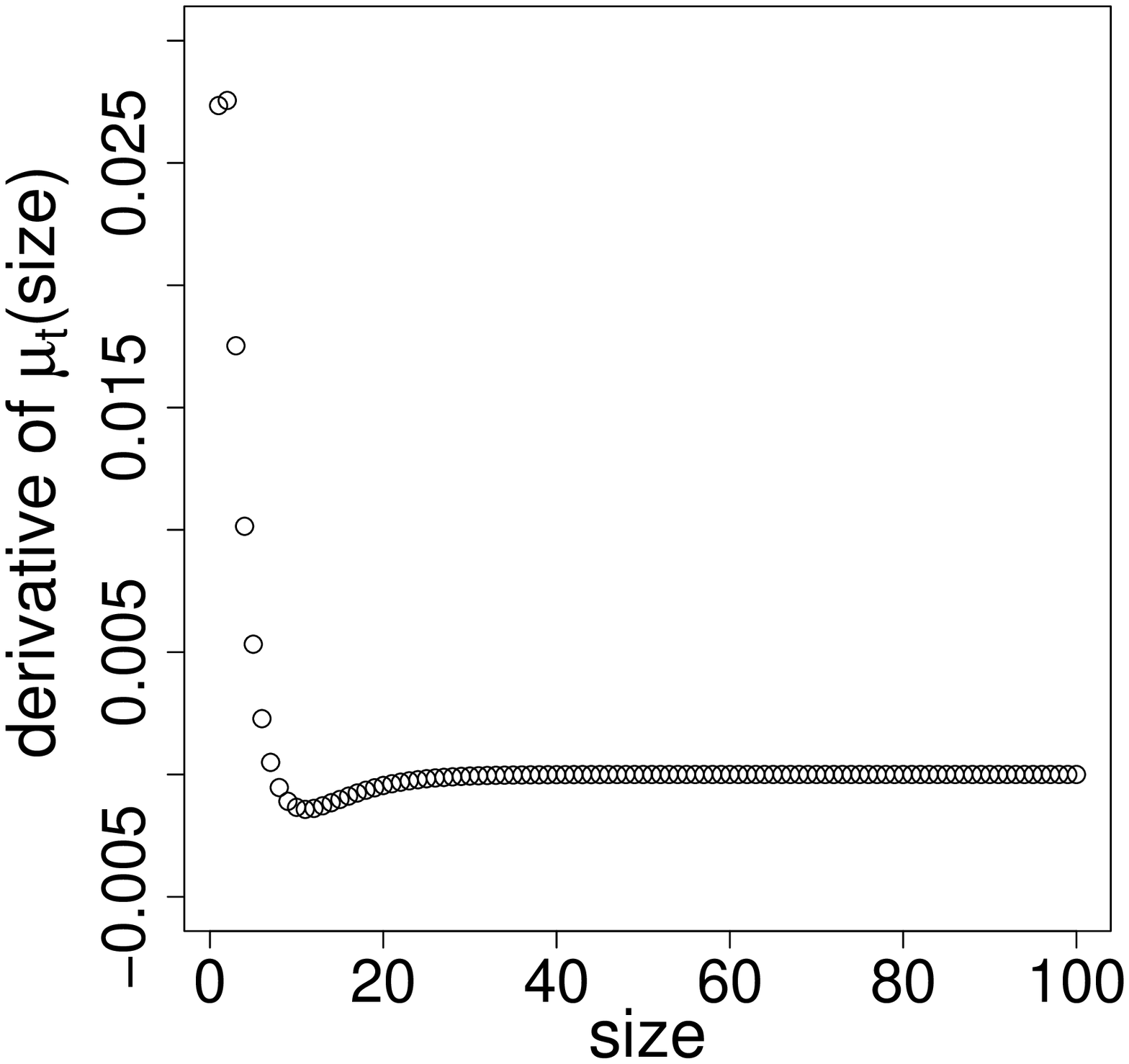}}
    }
    \subfigure[][$t=4.0$]
    {
        \resizebox{!}{0.35\linewidth}{\includegraphics{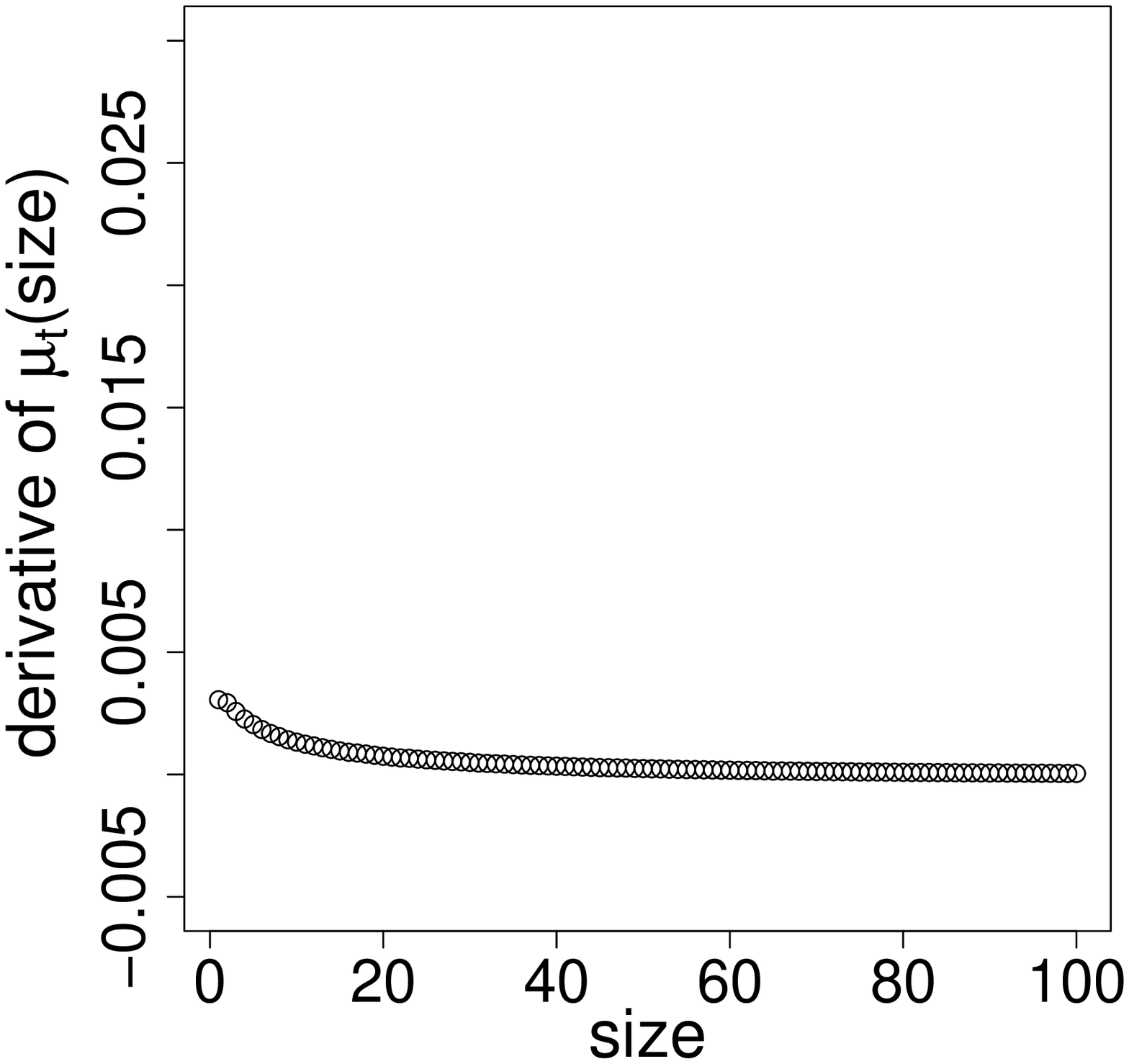}}
    }
\end{center}
    \caption{Derivative $\frac{\partial}{\partial \lambda} \mu_t^{\lambda}(x)$ versus particle size $x$, for Soot kernel using the Double Coupling algorithm, $\lambda=2.1$, $\epsilon=0.03$, $N=10^7$, $L=300$. Confidence intervals have been omitted since they are visually negligible.}
    \label{Soot:partsizeplot}
\end{figure}

\begin{figure}[!h]
\begin{center}
    \subfigure[][$t=0.5$]
    {
        \resizebox{!}{0.35\linewidth}{\includegraphics{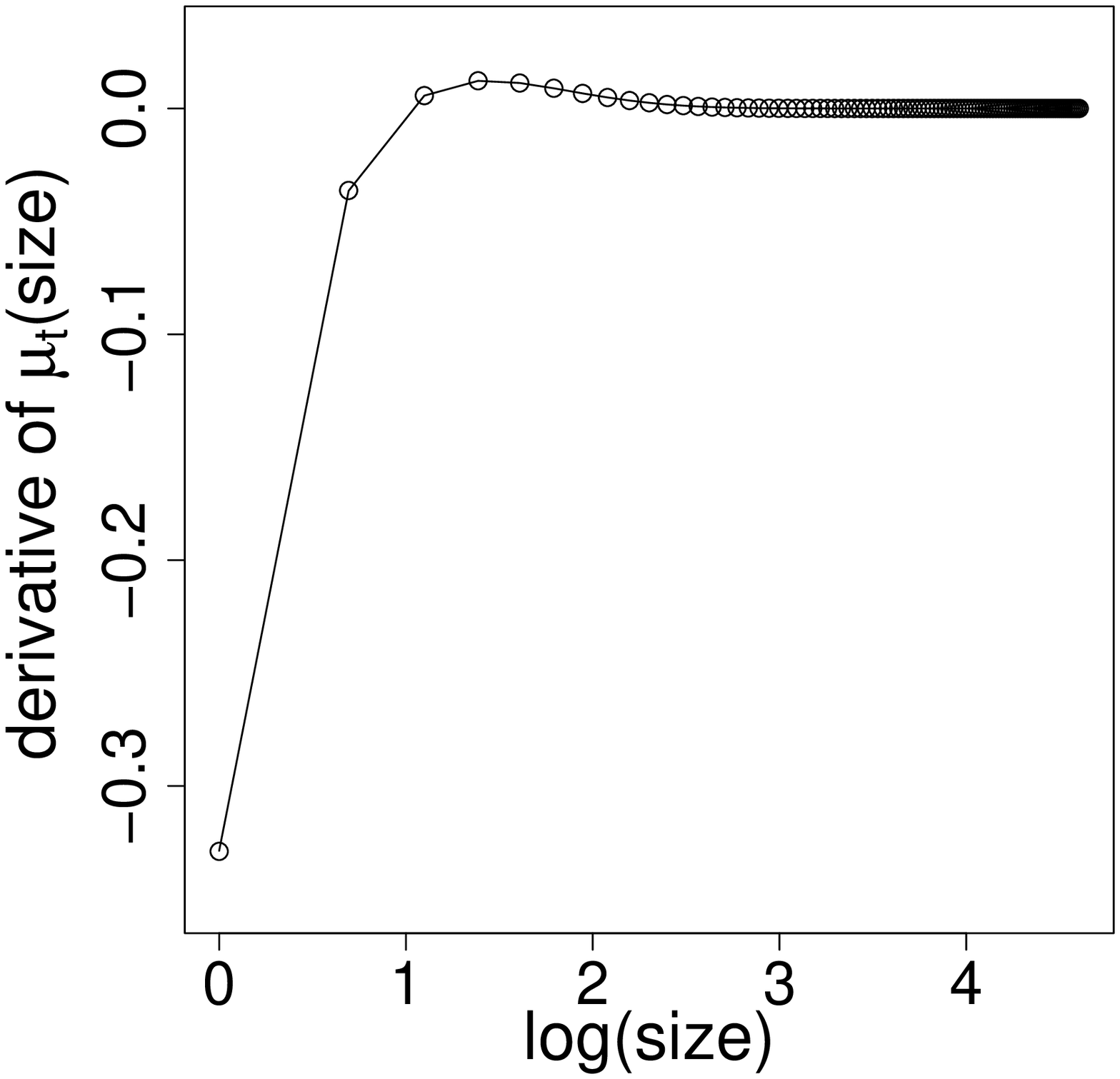}}
    }
    \subfigure[][$t=2.0$]
    {
        \resizebox{!}{0.35\linewidth}{\includegraphics{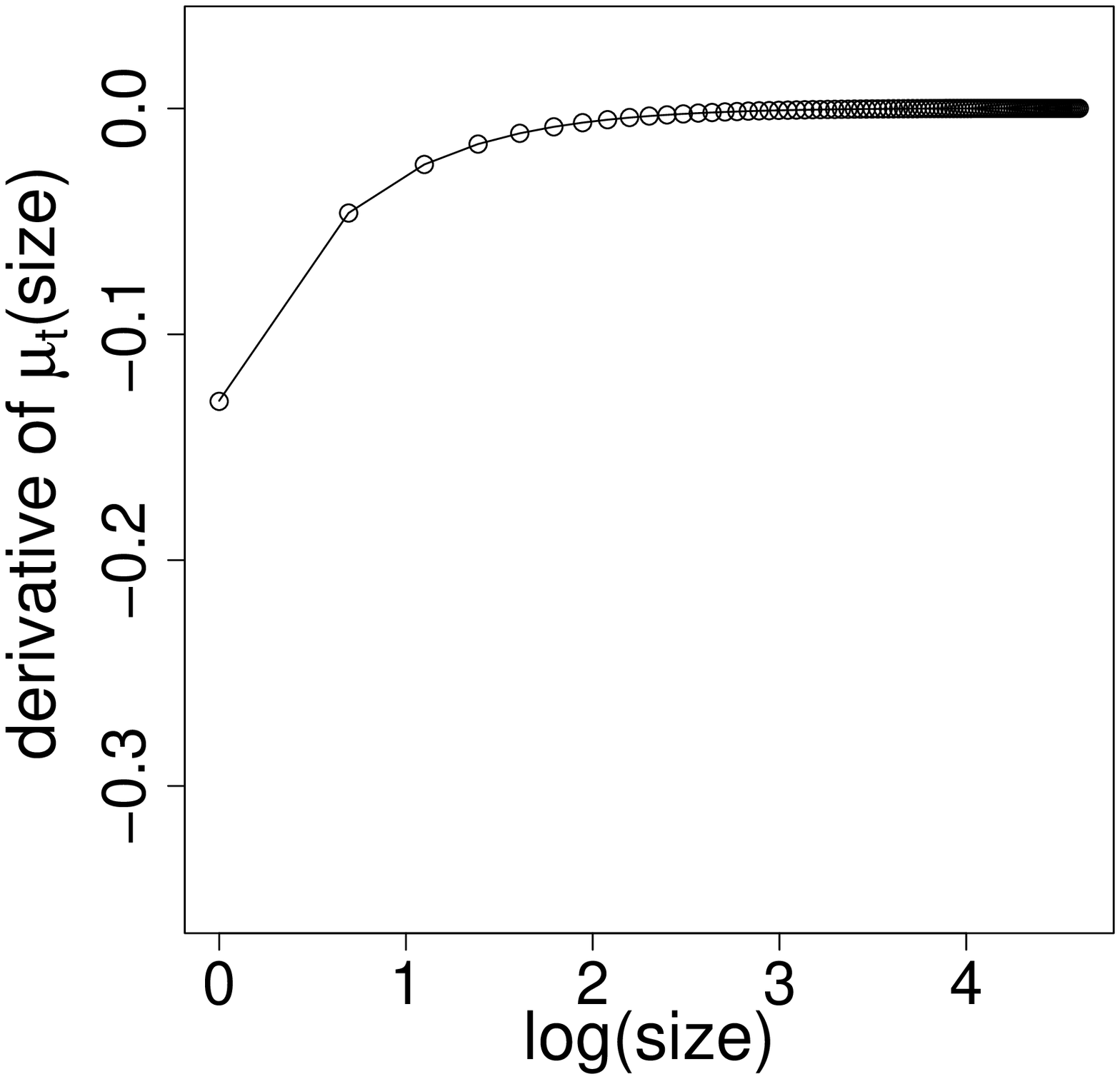}}
    }
\end{center}
    \caption{Derivative $\frac{\partial}{\partial \lambda} \mu_t^{\lambda}(x)$ versus particle size $x$, for additive kernel using the Double Coupling algorithm, $\lambda=1$, $\epsilon=0.06$, $N=10^6$, $L=300$. The line joins the points of the derivative of the analytic solution of eq.~\eqref{SmoEq} with monodisperse initial conditions, as given in \cite{Wagner3}. Confidence intervals have been omitted since they are visually negligible.}
    \label{Additive:partsizeplot}
\end{figure}

Figure \ref{Soot:partsizeplot} shows what the derivative of the parametric solution of $\mu_t^{\lambda}(x)$ looks like for the Soot kernel for two different evolution times $t$. Figure \ref{Additive:partsizeplot} shows similar quantities, but for the Additive kernel; it is in good agreement with the analytic solution given in \cite{Wagner3}.

\subsection{Convergence study}
There are two sources of systematic error in using the central difference estimator --- one due to using a non-zero value of $\epsilon$ and the other due to assuming a finite particle system. It is the latter we investigate --- here we estimate the order of convergence of the systematic error as $N$ varies. The value of $\epsilon$ will be fixed here.\\

We define the \textbf{systematic error due to $N$} as the difference between the expected central difference (for finite particle number $N$) and the analytic central difference.
\begin{equation}
e_{\textrm{sys}}(N;\lambda,\epsilon,t) = \mathbbm{E}\overline{F}^{(\lambda,\epsilon)}(N;t) - f^{(\lambda,\epsilon)}(t)
\end{equation}
However, as we nearly always do not know the analytic central difference $f^{(\lambda,\epsilon)}$, we estimate it using
\begin{equation}
\overline{F}^{(\lambda,\epsilon)}(N_{\textrm{large}};t)
\end{equation}
for very large $N_{\textrm{large}}$. Also, $\mathbbm{E}\overline{F}^{(\lambda,\epsilon)}(N;t)$ is estimated by $\overline{F}^{(\lambda,\epsilon)}(N;t)$. Now we set test function $f(y) = \mathbbm{1}_{\{y = i\}}$ again, to ensure that $\overline{F}^{(\lambda,\epsilon)}(N;t)$ is an estimate of the number density for particle size $i$. Therefore we rename $\overline{F}^{(\lambda,\epsilon)}(N;t)$ as $\overline{F}^{(\lambda,\epsilon,i)}(N;t)$ for number density estimate at particle size $i \in \mathbbm{N}$; we adopt analogous notation $f^{(\lambda,\epsilon,i)}$ and $e_{\textrm{sys}}(N;\lambda,\epsilon,t,i)$ for this particular choice of test function. Our metric for considering convergence in $N$ is simply the absolute estimated systematic error, summed over chosen evolution times $(t_k)_{k=1}^{T}$\footnote{We take the $t_k$ to be $(0.5,1.0,\ldots,7.0)$} and summed over particle sizes $i$:
\begin{equation}
c_{\textrm{tot}} = \sum_{k=1}^T \sum_{i \in \mathbbm{N}} |e_{\textrm{sys}}(N;\lambda,\epsilon,t_k,i)|
\end{equation}

\begin{figure}
\begin{center}
    \subfigure[][Additive kernel, $\lambda = 1$, $\epsilon=0.06$ ]
    {
        \resizebox{!}{0.35\linewidth}{\includegraphics{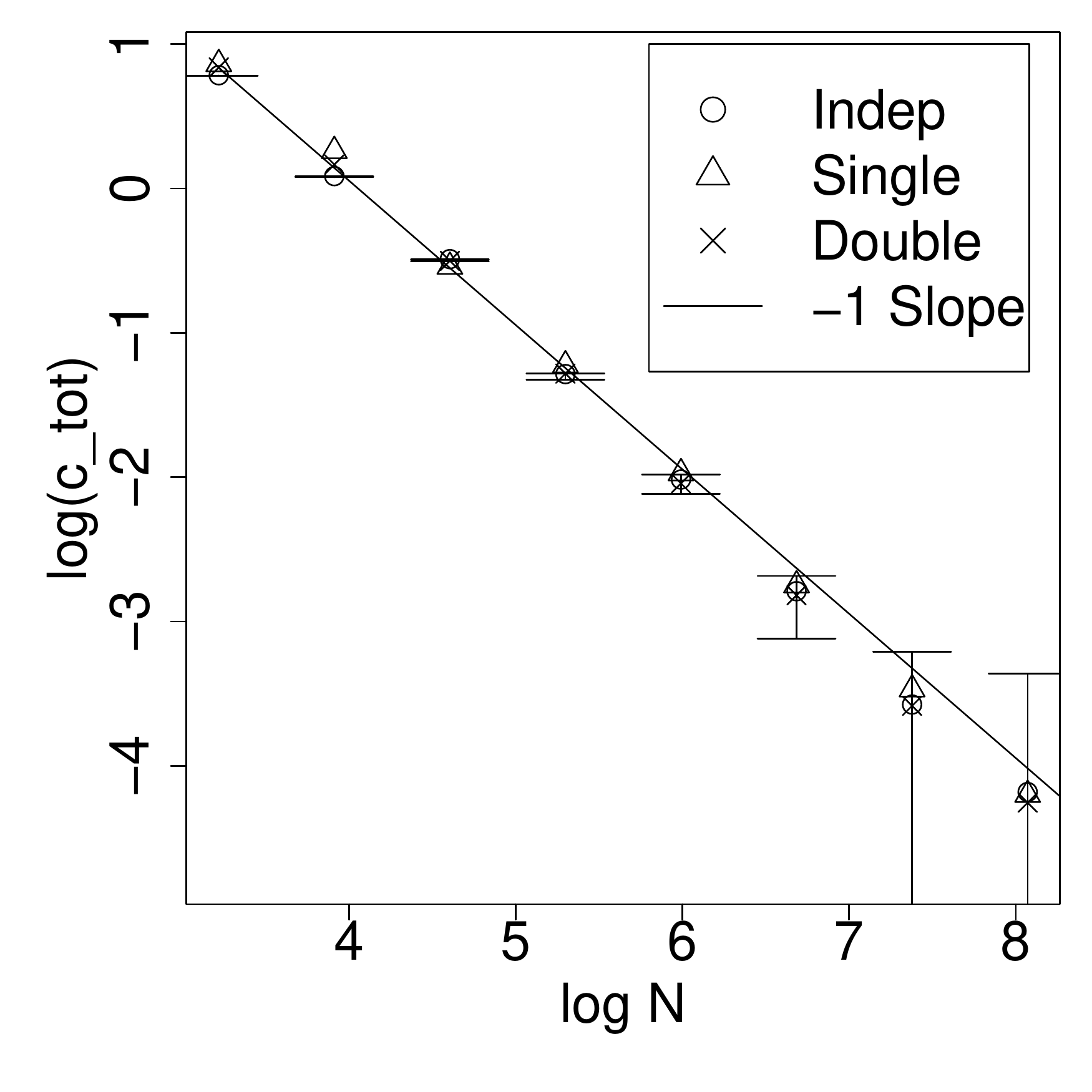}}
    }
    \subfigure[][Soot kernel, $\lambda = 2.1$, $\epsilon=0.03$]
    {
        \resizebox{!}{0.35\linewidth}{\includegraphics{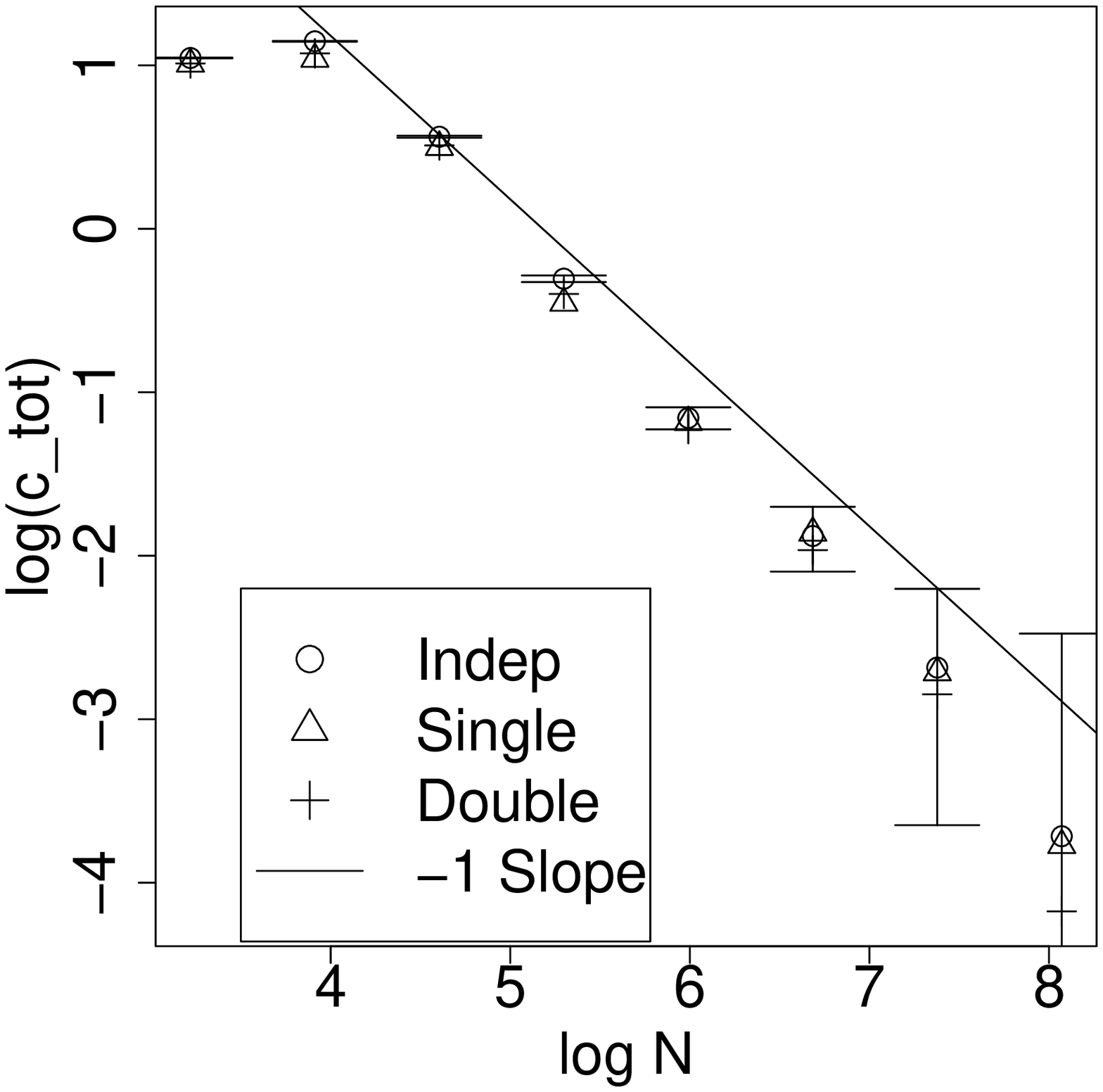}}
    }
\end{center}
    \caption{$\log c_{\textrm{tot}}$ versus $\log N$, $N=25 \times 2^i$ for $i=0,\ldots,7$, $N \times L=10^8$ and $0 \leq t \leq 3$. Confidence intervals given are for the Independent case only.}
    \label{convergenceplots}
\end{figure}
\noindent Figure \ref{convergenceplots} shows what we expect --- the $c_{\textrm{tot}} \sim \frac{1}{N}$. We obtain similar plots for different values of $\epsilon$.


\subsection{Statistical error}
\label{stat_err_section}

\begin{figure}[!h]
\begin{center}
    \subfigure[][$t=1.0$]
    {
        \resizebox{!}{0.35\linewidth}{\includegraphics{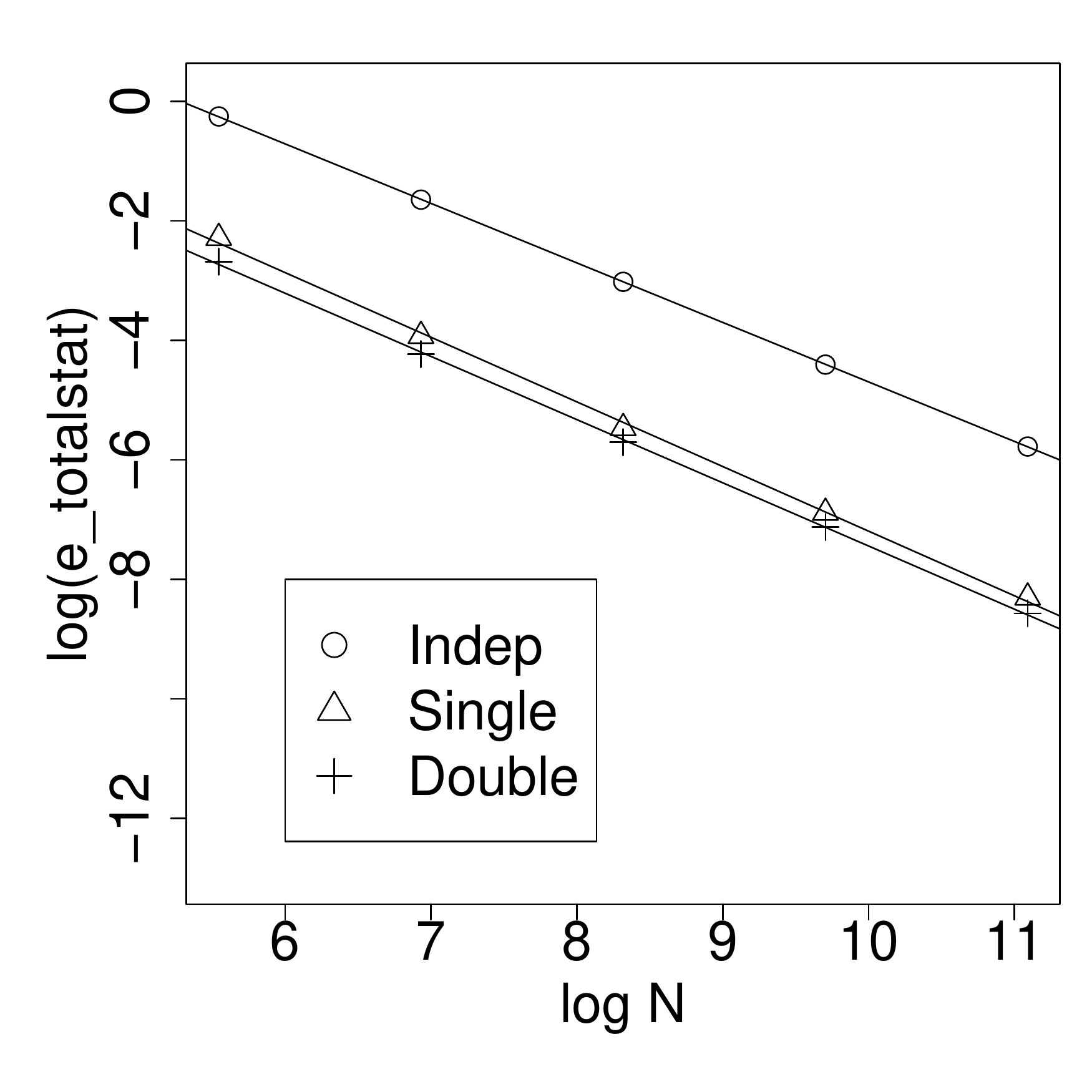}}
    }
    \subfigure[][$t=2.0$]
    {
        \resizebox{!}{0.35\linewidth}{\includegraphics{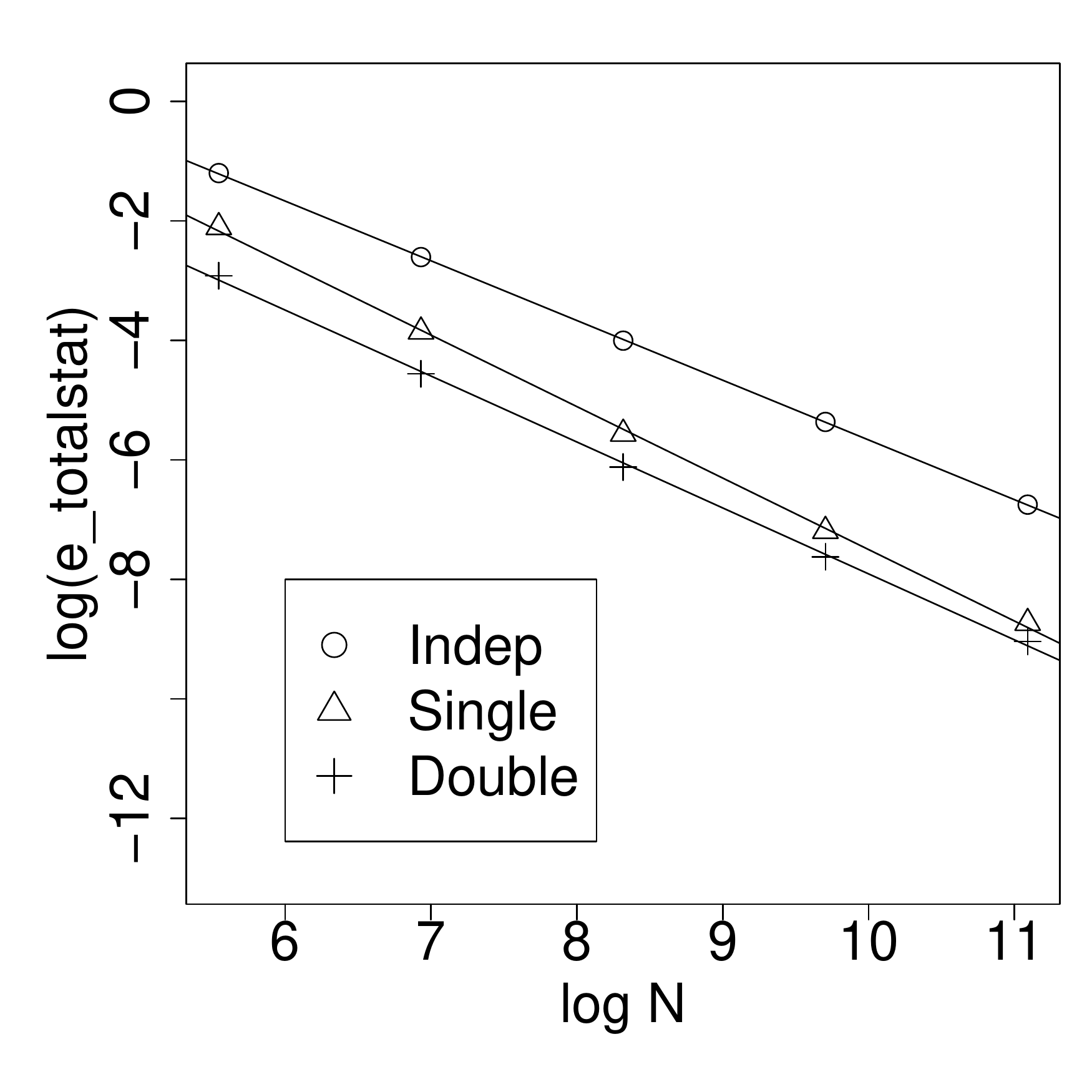}}
    }\\
    \subfigure[][$t=3.0$]
    {
        \resizebox{!}{0.35\linewidth}{\includegraphics{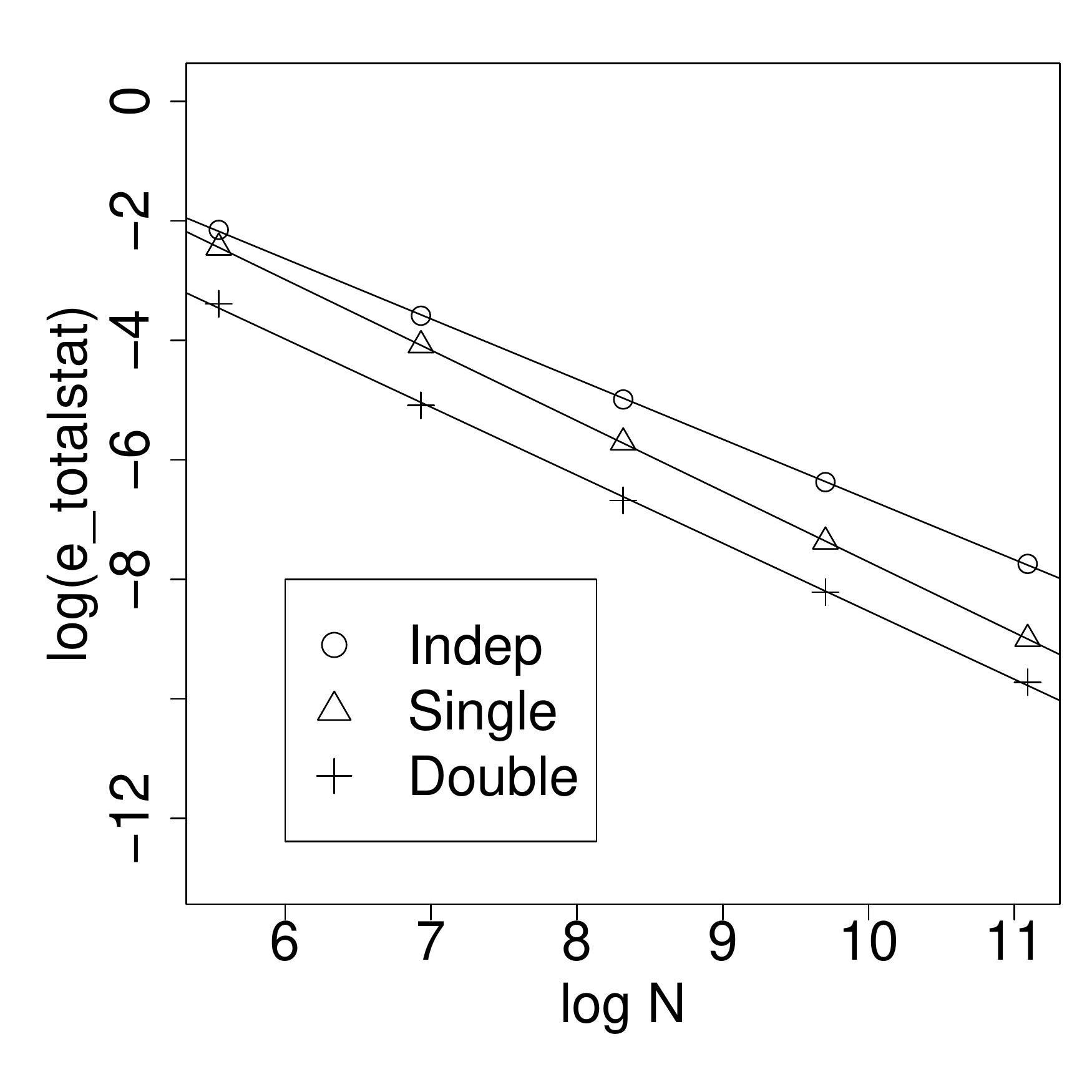}}
    }
    \subfigure[][$t=4.0$]
    {
        \resizebox{!}{0.35\linewidth}{\includegraphics{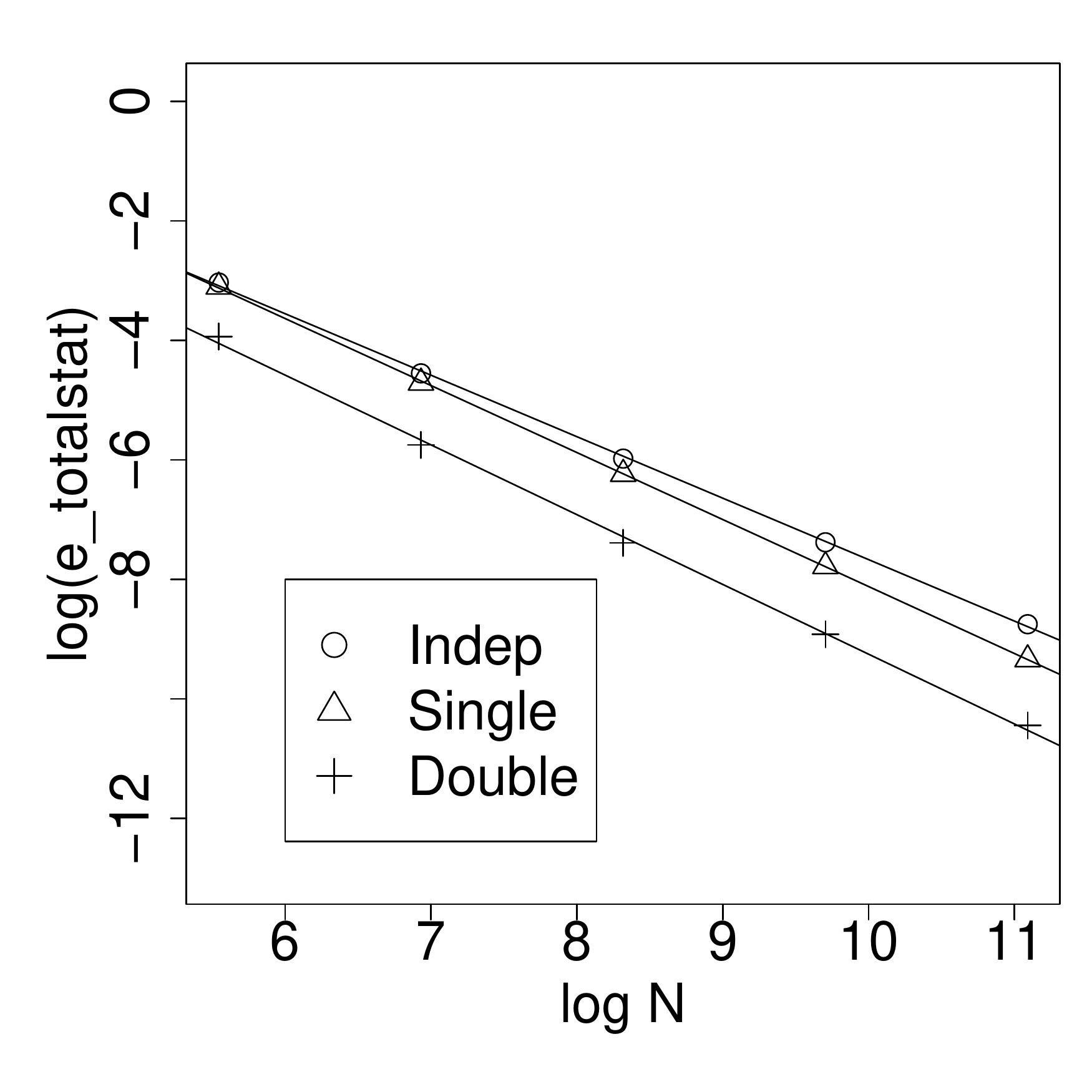}}
    }\\
\end{center}
    \caption{$\log e_{\rm{totalstat}}$ versus $\log N$ over different values of $t$ for all the algorithms for the Additive kernel, where $N=2^{10},2^{12},2^{14},2^{16},2^{18}$,$\, N \times L=2^{26}$ and $\epsilon=0.06$.}
    \label{Additive:staterr_study}
\end{figure}

This quantity is defined as
\begin{equation}
e_{\textrm{stat}}(N;t,\lambda,\epsilon) = \overline{F}^{(\lambda,\epsilon)}(N;t) - \mathbbm{E} \overline{F}^{(\lambda,\epsilon)}(N;t).
\end{equation}
This is a signed measure which associates to each $i\geq 1$ the number $e_{\textrm{stat}}(N;t,\lambda,\epsilon,i)$. We consider in this paragraph how it behaves according to the different algorithms and kernels. The variance of each estimator $F_l$ can be estimated for each $i\geq 1$ by
\begin{equation}
v_F^{(i)} := \frac{1}{L-1} \sum_{l=1}^{L} (F_l(i) - \overline{F}(i))^2.
\end{equation}
This implies that the asymptotic $100(1-\alpha) \%$ confidence interval for $\mathbbm{E} \overline{F}^{(\lambda,\epsilon)}(N;t;i)$ is:
\begin{equation}
\left[ \overline{F}^{(\lambda,\epsilon)}(N;t;i) - z_{\alpha / 2} \sqrt{\frac{v_F^{(i)}}{L}} \: , \: \overline{F}^{(\lambda,\epsilon)}(N;t;i) + z_{\alpha / 2}\sqrt{\frac{v_F^{(i)}}{L}} \right]
\end{equation}
where $z_{\alpha /2}$ is the upper $\alpha / 2$ point of the standard normal distribution. Hence,
\begin{equation}
\mathbbm{P} \Bigl( \big|e_{\textrm{stat}}(N;t,\lambda,\epsilon)\big| \leq z_{\alpha / 2}\sqrt{\frac{v_F^{(i)}}{L}} \Bigr) \approx 1-\alpha.
\end{equation}
For this paper, we set $\alpha = 0.05$ \ie we consider $0.95\%$ confidence intervals. Also, consider the sum $e_{\textrm{totalstat}}$ of the single-simulation variances over the particle sizes:
\begin{equation}
e_{\textrm{totalstat}}:=\sum_{i\geq 1} v^{(i)}_{F}
\end{equation}
We wish to see how this quantity behaves with $N$. Figure \ref{Additive:staterr_study} demonstrates that for all three algorithms, the total variance $e_{\textrm{totalstat}}$ behaves as $\frac{1}{N}$ since the slopes of the fitted lines are approximately $-1$. More importantly, the intercept for the Double algorithm is lower than the Single case, and much lower than that of the Independent case indicating that $e_{\textrm{totalstat}}$ is much smaller for the Double and Single cases than for the Independent---$e_{\textrm{totalstat}}$ for the Independent case at one point is approximately 20 times larger than that for the Double case at $t = 1.0$. One interesting observation is that the difference in the intercepts (of the fitted lines) between the Double and Indep decreases, showing that the benefits of smaller statistical error in the Double case become less pronounced as $t$ increases. The same happens for the Single algorithm, but at a faster rate, showing that the $\Xminus$ and $\Xplus$ systems diverge from each other, but faster for the Single algorithm than for the Double algorithm. 

\subsection{Efficiency}
\begin{figure}[!h]
    \subfigure[][$N=10^2$, $\epsilon=0.01$]
    {
        \includegraphics[width=0.3\textwidth]{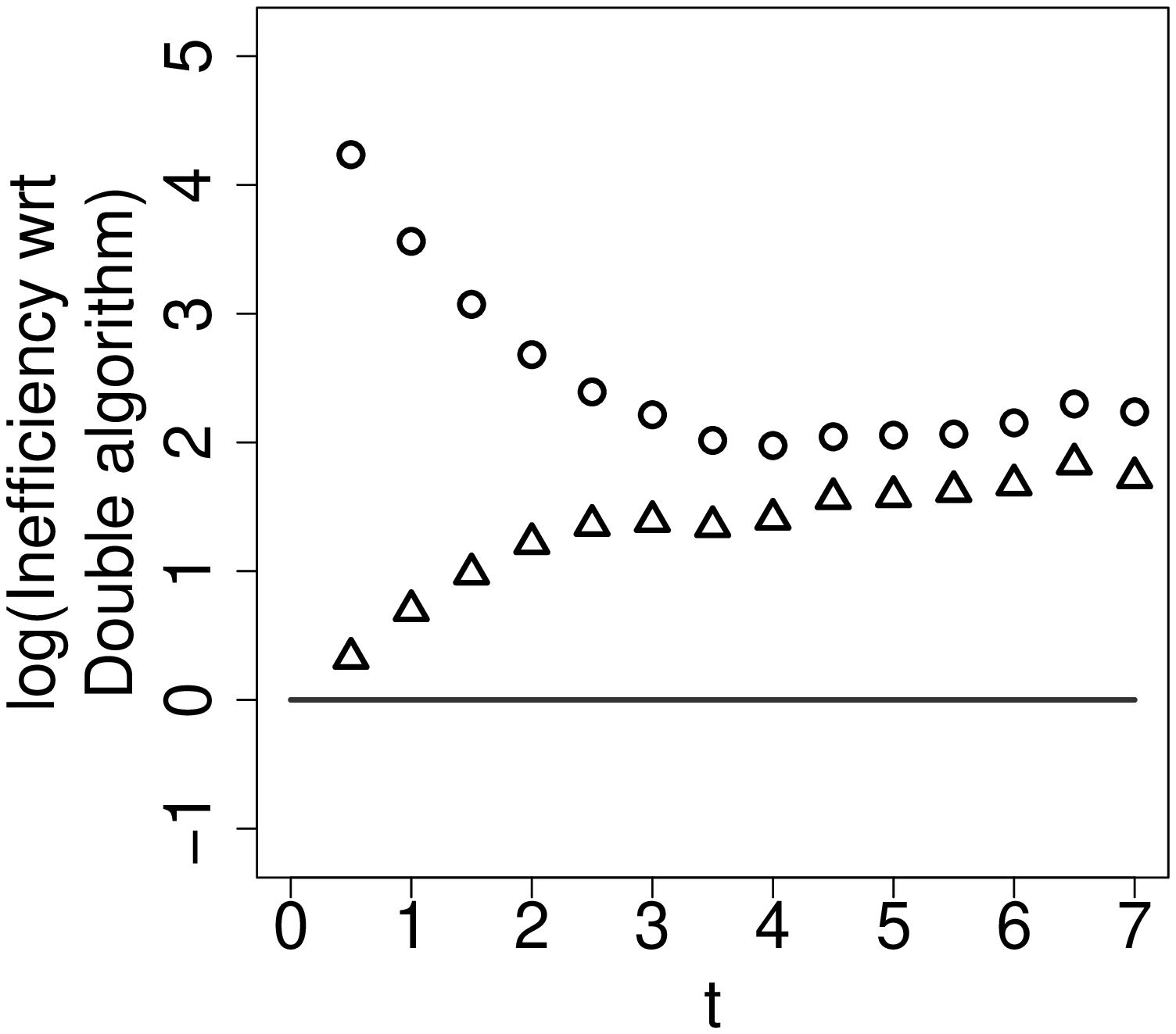}\\
    }
    \subfigure[][$N=10^2$, $\epsilon=0.05$]
    {
        \includegraphics[width=0.3\textwidth]{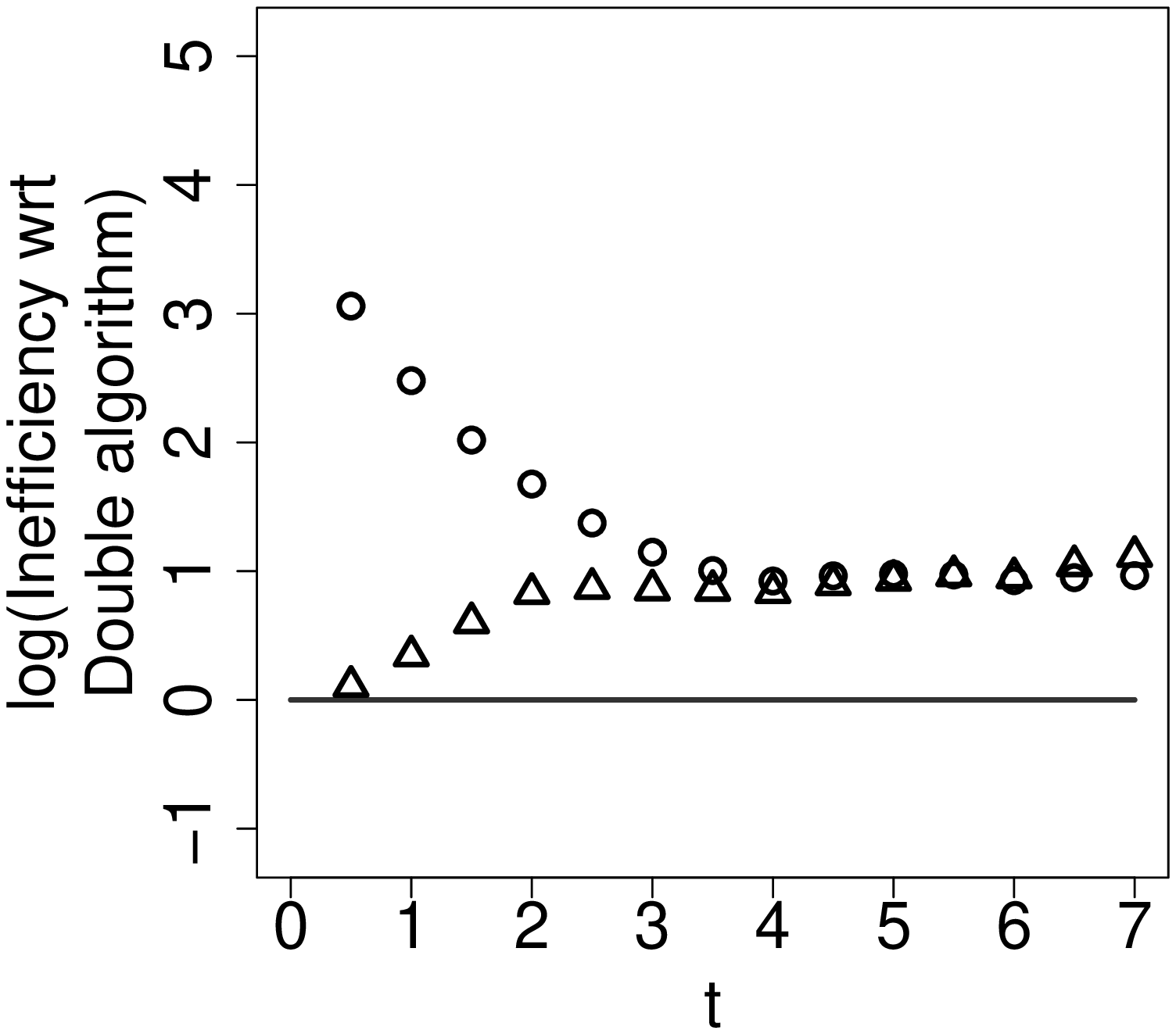}\\
    }
    \subfigure[][$N=10^2$, $\epsilon=0.20$]
    {
        \includegraphics[width=0.3\textwidth]{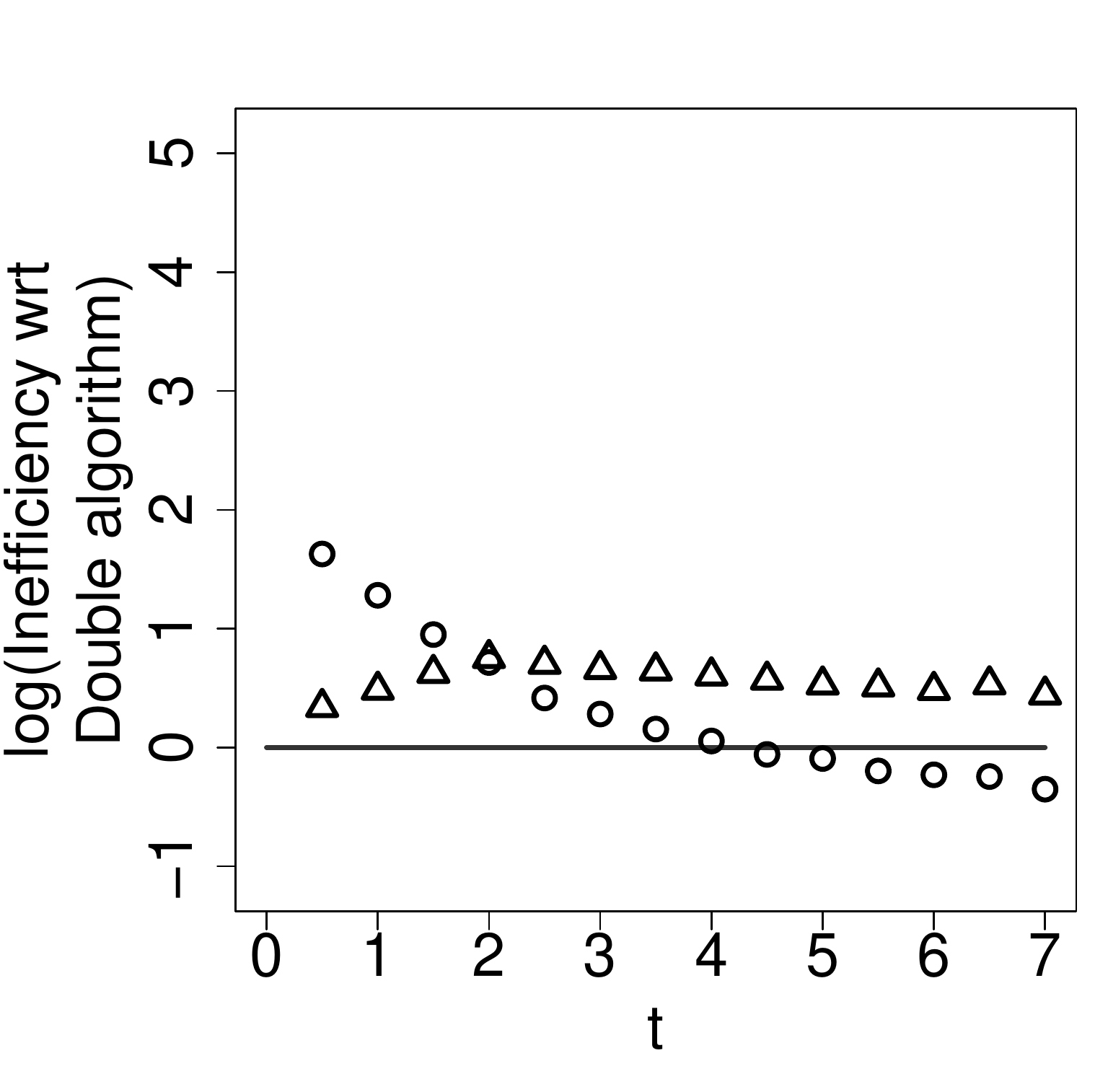}\\
    }\\
    \subfigure[][$N=10^5$, $\epsilon=0.01$]
    {
        \includegraphics[width=0.3\textwidth]{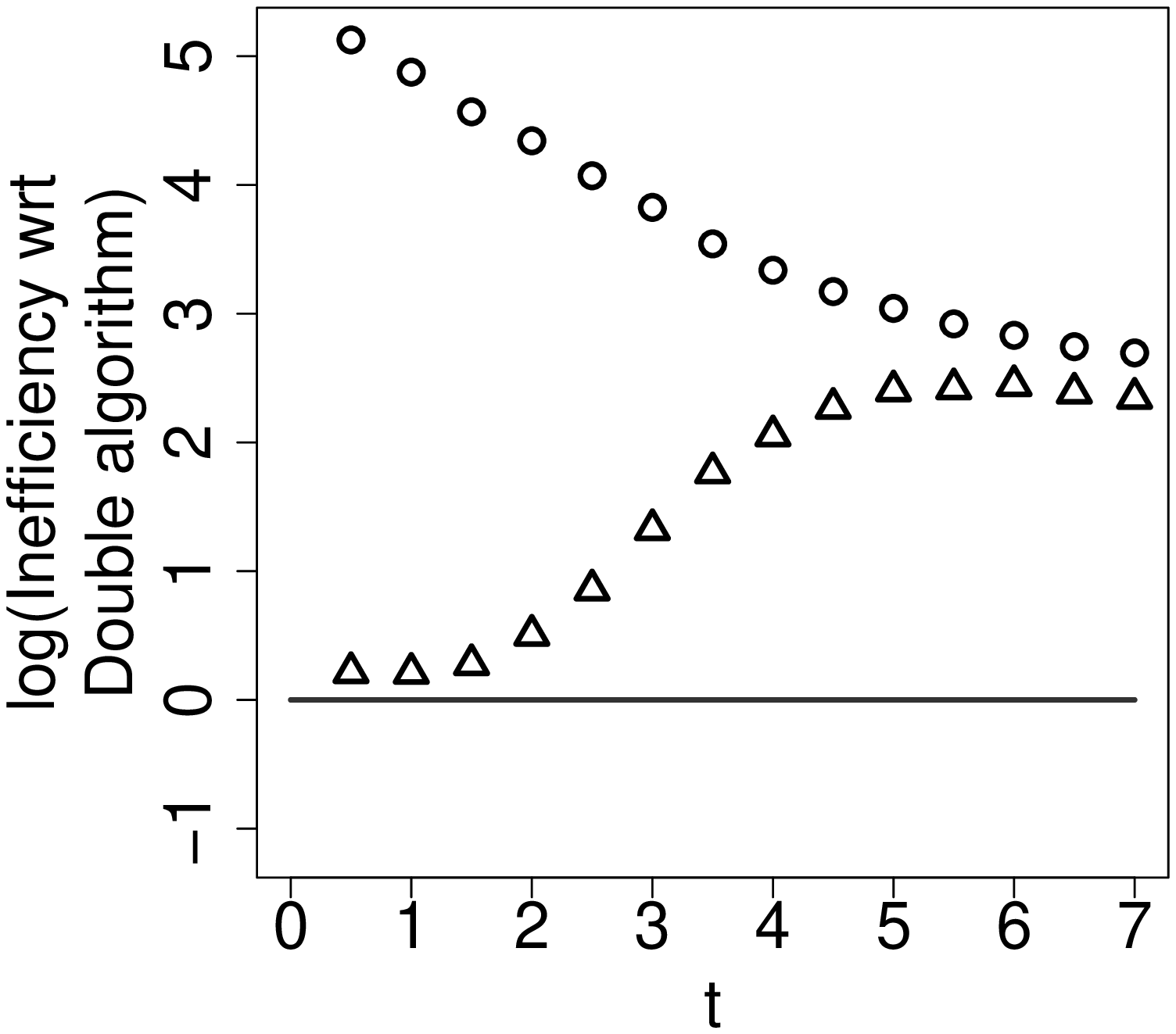}\\
    }
    \subfigure[][$N=10^5$, $\epsilon=0.05$]
    {
        \includegraphics[width=0.3\textwidth]{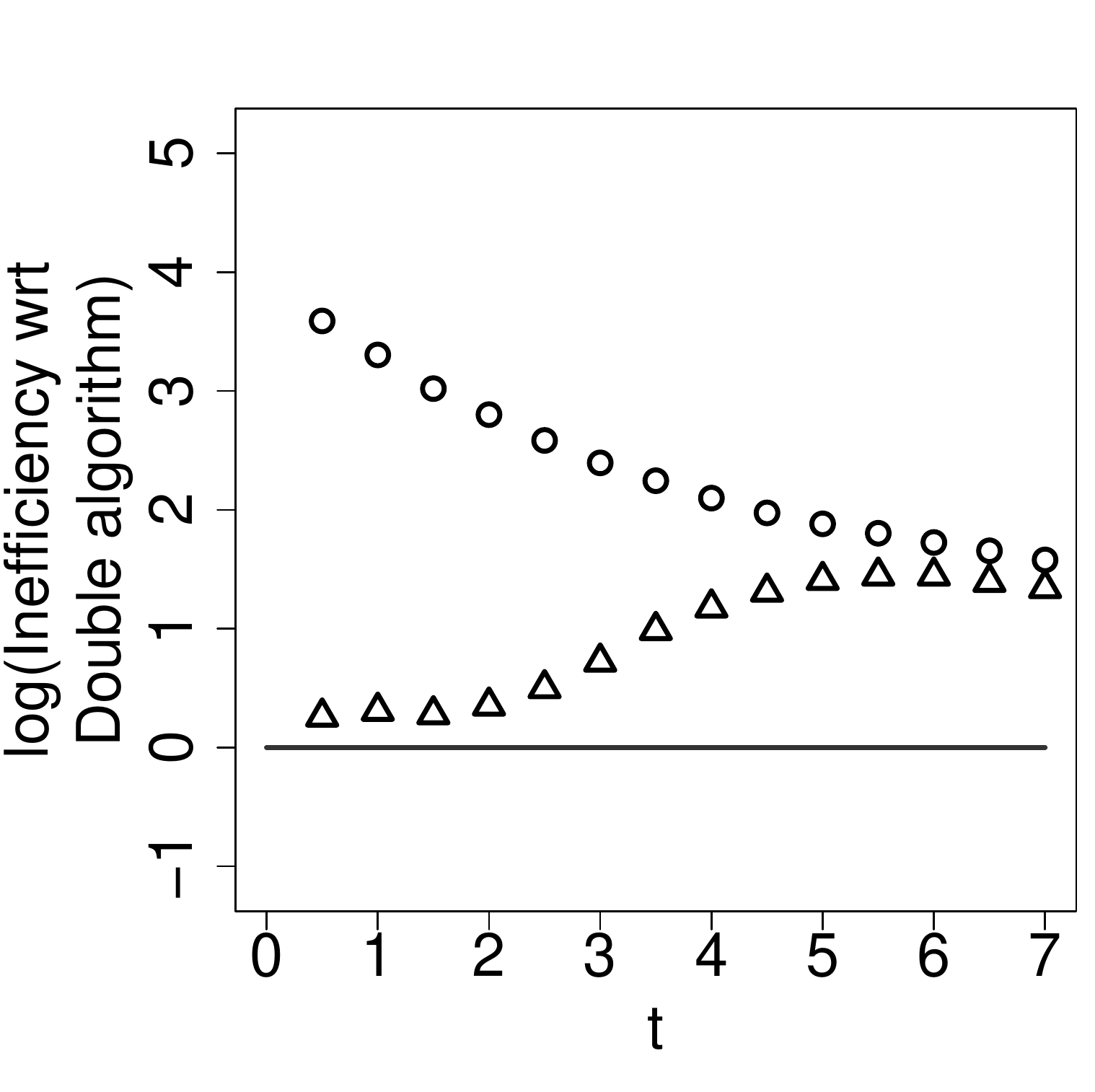}\\
    }
    \subfigure[][$N=10^5$, $\epsilon=0.20$]
    {
        \includegraphics[width=0.3\textwidth]{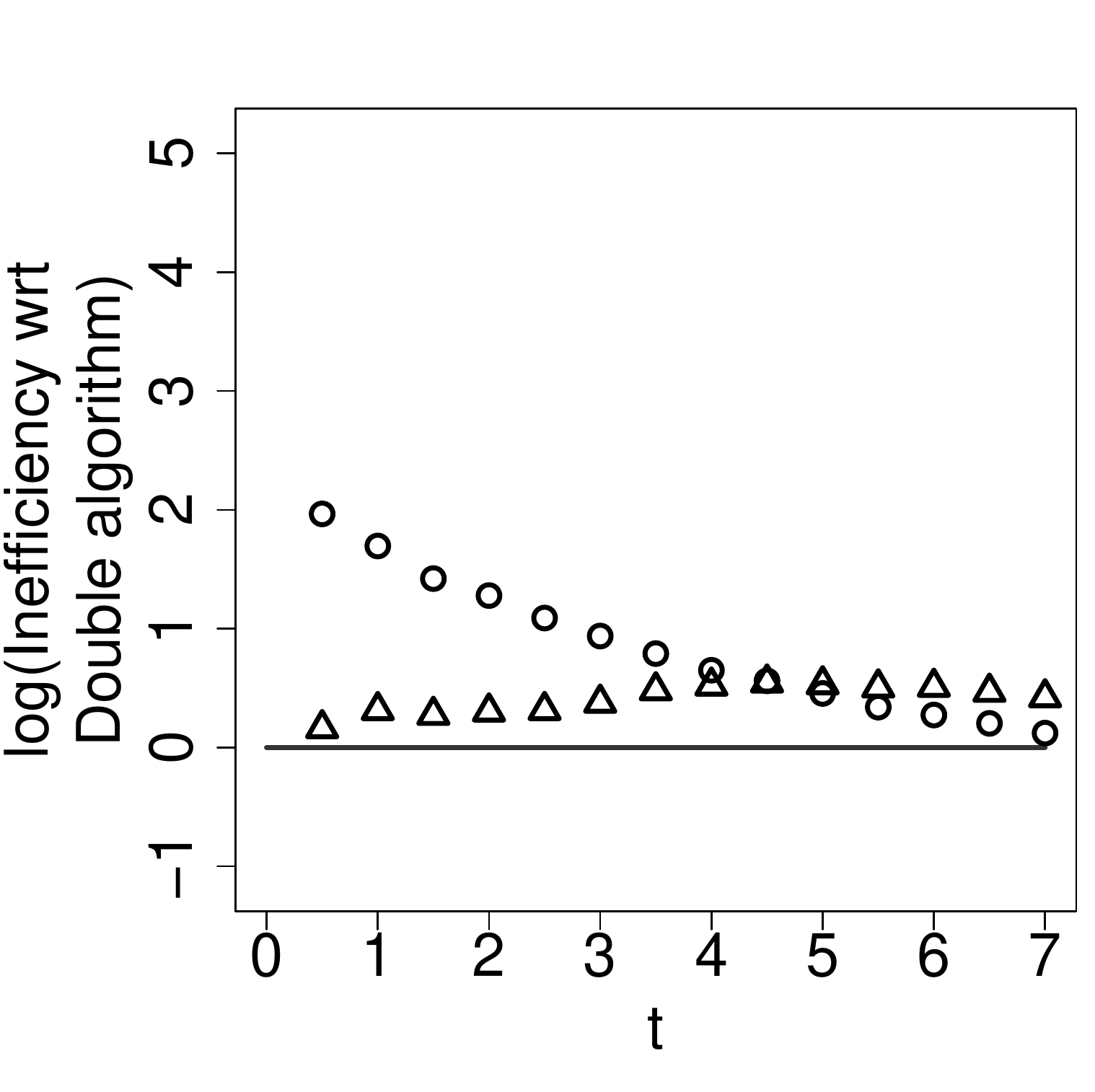}\\
    }
    \caption{Additive kernel --- $\textrm{Inefficiency relative to Double}$ for Indep and Single algorithms, as a function of $t$ for different values of $\epsilon$ and $N$. The Independent algorithm is represented by circles, Single by triangles and the Double threshold by the horizontal line.}
    \label{Eff_additive_graphs}
\end{figure}

In this subsection, we will discuss which algorithm is `best'. A quantification of this quality needs to be defined which will take into account the accuracy and the run time of each algorithm. One such measure can be described as the run time needed to achieve a certain fixed statistical error (assuming $N$ and $\epsilon$ are fixed).
Before defining it, we shall
\begin{romannum}
\item Let $t_{\text{run}}(t)$ be the CPU time actually taken to perform $L_{\text{run}}$ simulations.
\item Set $e_{\textrm{fixed}}=\frac{e_{\textrm{totalstat}}}{L}$; this ``total standard error'' will be artificially \textbf{fixed}.
\item Let $L_{\textrm{est}}(t)$ be the estimated number of simulations required to acquire the fixed $e_{\textrm{fixed}}(t)$ value and $t_{\textrm{est}}(t)$ be the CPU time required to perform $L_{\textrm{est}}(t)$ runs.
\end{romannum}
The condition $e_{\textrm{fixed}}(t) = \sum_{i\geq 1} \frac{v^{(i)}_{F}}{L_{\textrm{est}}(t)}$ implies that $L_{\textrm{est}}(t) = \frac{ \sum_{i\geq 1} v^{(i)}_{F}}  {e_{\textrm{fixed}}}$; we then have that $t_{\textrm{est}}(t)=\frac{t_{\text{run}}(t)}{L_{\text{run}}(t)}L_{\textrm{est}}(t)$. We shall use the following quantity to compare different algorithms.
\begin{equation}
\textrm{Inefficiency}^{\textrm{algorithm}}
:= \frac{t^{\textrm{algorithm}}_{\textrm{est}}(t)}{t^{\textrm{Double}}_{\textrm{est}}(t)};
\end{equation}
we call it the \textbf{Inefficiency with respect to the Double Coupling algorithm}. Hence, if an algorithm has an Inefficiency of more than unity, the algorithm does not perform as well as the Double Coupling algorithm. Figure \ref{Eff_additive_graphs} plots these inefficiencies.
One can see that the Independent algorithm has large inefficiencies for small $t$---this is due to the vastly smaller statistical errors of the Double and Single Coupling algorithms, as well as all three central difference algorithms taking comparable times to run. In fact, the Double and Single algorithms are generally quicker for smaller $\epsilon$ since the Independent algorithm requires two simulations to generate a derivative estimate, as well as the fact that the Single and Double algorithms have fewer $\ominus$ and $\oplus$ particles to deal with. The inefficiencies of the Single algorithm lie between 1.0 and 2.0 (note that Figure \ref{Eff_additive_graphs} uses log scales) indicating that the Double algorithm has a significant improvement over the Single in terms of accuracy. Also, the inefficiencies decrease with $\epsilon$ since larger $\epsilon$ implies that $\oplus$ and $\ominus$ particles are created, thus meaning that both coupling algorithms have larger CPU run times (the run times increase almost linearly with $\epsilon$ for the Single and Double, whereas they are almost constant with respect to $\epsilon$ for the Independent algorithm). Also, the ratio of variances for the Independent and Single algorithms relative to the Double algorithm decreases for larger $\epsilon$, due to the increasing similarity between the Independent algorithm and the coupling algorithms, thus decreasing the inefficiencies. Note however that the case $\epsilon=0.2$ is unlikely to be computationally useful as it amounts to a 20\% change in the parameter, which enters multiplicatively into the kernel.

Furthermore, one notices that the larger value of $N=10^5$
results in larger inefficiencies for the Independent and Single algorithms for the Additive kernel. This appears to be because as $N$ increases, the capacity for cancellations is larger, implying that the Double is more accurate (and faster) than one
expects. Also, as $t$ increases, we find that the number of particles for all algorithms decreases dramatically, and so there appears to be little difference in accuracy between all three algorithms for larger $t$. This is to be expected since $\mu_t^{+}$ and $\mu_t^{-}$ will become increasing dissimilar, thus indicating smaller covariances, and so the variances for the Double and Single Coupling algorithms will be similar to those for the Independent algorithm, this loss of efficiency of the coupling algorithms is ultimately unavoidable in this class of algorithms, and one can only hope to minimise this decrease.

It is important to realise that this analysis does not take into account the systematic error due to $\epsilon$ or $N$ since the Inefficiency metric only uses estimated variances. A related problem with the analysis is that the number of particles for $t \in [3.5,7.0]$ becomes quite small\footnote{about $O(1)-O(10)$ if $N=10^2$ but $O(100)-O(1000)$ for $N=10^6$.}, and therefore the systematic errors and estimated variances are not very reliable.

\section{Conclusions}
In this paper, two new stochastic algorithms were described which solve for parametric derivatives of the solution to the discrete Smoluchowski's coagulation equation. These algorithms consider two Marcus-Lushnikov processes which are coupled together in order to reduce the difference in their trajectories. The hope was that this would significantly reduce the variance of the central difference estimators of the parametric derivatives. In the numerical results section, we first validated the fact that the order of convergence for these algorithms is indeed $O(1/N)$. Furthermore, it was shown from the statistical error plots that the accuracy is order of magnitudes better than that of the worst case (the Independent algorithm), at least for larger $N$ and smaller $\epsilon$. Subsequently, we considered a method of comparing the algorithms which considers both the variances of the derivative estimators as well as the CPU run times. It was shown that the Double algorithm is mostly more `efficient' than Single over variations in $\epsilon$ and $t$, whilst being significantly more `efficient' than the Independent algorithm for small $t$, large $N$ and small $\epsilon$, though some of this advantage is lost for larger $t$ and $\epsilon$.
%

\smallskip

%
%

%
%

\end{document}